\documentclass[11pt]{article}
\usepackage[utf8]{inputenc}

\usepackage{amsmath,amsthm,graphicx,amssymb,float,authblk,color,subcaption,appendix,bm}
\usepackage{algorithm,algpseudocode}
\usepackage{url,hyperref}
\usepackage[round]{natbib}
\usepackage[sc]{mathpazo}
\usepackage[gen]{eurosym}
\usepackage[capitalise]{cleveref} 
\usepackage[dvipsnames]{xcolor}
\usepackage{placeins}

\usepackage{setspace}
\usepackage{booktabs}
\usepackage{threeparttable}
\newtheorem{proposition}{Proposition}
\newtheorem{theorem}{Theorem}
\newtheorem{lemma}{Lemma}
\newtheorem{definition}{Definition}
\newtheorem{corollary}{Corollary}
\newtheorem{assumption}{Assumption}
\newtheorem{remark}{Remark}

\newcommand{\setS}{\mathcal{S}}
\newcommand{\setA}{\mathcal{A}}
\newcommand{\setN}{\mathcal{N}}

\newcommand{\setG}{\mathcal{G}}
\newcommand{\setL}{\mathcal{L}}

\newcommand{\bbR}{\mathbb{R}}
\newcommand{\bbE}{\mathbb{E}}

\newcommand{\vecs}{\mathbf{s}}
\newcommand{\veca}{\mathbf{a}}

\newcommand{\vecpi}{\bm{\pi}}

\newcommand{\setT}{\mathcal{T}}

\title{Markov Population Potential Game: An Alternative Framework for Markovian Traffic Assignment}
\author[1]{Xuhang Liu}
\author[1]{Kenan Zhang\footnote{Corresponding author \url{kenan.zhang@epfl.ch}}}
\author[2]{Rui Yao}
\affil[1]{School of Architecture, Civil and Environmental Engineering,
École polytechnique fédérale de Lausanne, Lausanne, Switzerland}
\affil[2]{Department of Civil and Environmental Engineering, Technion – Israel Institute of Technology, Haifa, Israel}

\date{}

\begin{document}
\maketitle

\begin{abstract}
Understanding the equilibrium behavior of large populations of interacting agents in transportation systems is challenging due to the high-dimensional state space and the coupled decision-making among agents. This study proposes the Markov population potential game (MPPG) that integrates Markov games, potential games, and population games as a unified framework for analyzing such interactions. 
The framework extends the classical Markov potential game to the population level and establishes theoretical conditions on stage rewards and state transitions under which the Markov population game admits a Beckmann-type potential function. Accordingly, every KKT point of the associated optimization problem yields an MPPG equilibrium.
We further demonstrate that if only the reward condition holds, each KKT point of the optimization problem corresponds to a support-restricted $\varepsilon$-equilibrium, 
meaning the Q-values of any two active actions differ by at most $\varepsilon$, thereby significantly extending the applicability of MPPG. We further establish both priori and 
posteriori upper bounds of $\varepsilon$ to efficiently evaluate the approximation error.
Inspired by the equivalent optimization problem, a projected-gradient (PG) policy iteration algorithm is developed. 
The theoretical bounds and PG algorithm are demonstrated through numerical experiments on a parametric Markov population game. 

\end{abstract}


\section{Introduction}

Modeling the behaviors of a large population of strategic agents (e.g., travelers and vehicles) is crucial for the design and management of transportation systems. While static traffic assignment models~\citep{wardrop1952road,beckmann1956studies,smith1979existence} are well-suited for long-term planning, emerging problems in transportation are characterized by more dynamic and complex interactions. 
Existing dynamic traffic assignment (DTA) models, although properly capturing traffic dynamics and propagation in the network, are mostly intractable~\citep{peeta2001foundations} and thus rely on simulation-based solution approaches~\citep{mahmassani1993network,peeta1995multiple,mahmassani2001dynamic}. 
The primary goal of this study is to develop an intermediary modeling approach that balances representational power and computational tractability by adopting Markovian dynamics. Specifically, we model the sequential decisions of agents in the transportation system and connect their collective behaviors to the overall system performance and dynamic evolutions. 
By assuming memoryless (Markovian) transitions for both individual and system states, the model largely reduces complexity while preserving key factors that drive dynamics in the transportation system. 

In the Markovian framework, each agent's decision-making process is defined as a Markov decision process (MDP). When interactions among agents influence each one's payoff and state transitions, their collective behaviors lead to a Markov game~\citep{shapley1953stochastic}.
A standard solution concept of Markov games is the Nash equilibrium (NE), where no agent can improve its payoff by unilaterally deviating from its strategy. However, designing a tractable algorithm to find the Nash equilibrium is in general challenging, particularly when the number of agents increases \citep{papadimitriou2007complexity}. 
Nevertheless, when agents are homogeneous and the population approaches infinity, the mean-field approximation can be employed to reduce the complexity of computing Nash equilibria ~\citep{tembine2009mean,gast2012mean,elokda2024dynamic}, replacing the interactions among agents by the
interaction between each agent and the mean field. The resulting games are widely known as mean field games~\citep{huang2006large, lasry2007mean}, or in the context of this paper, Markov population games with a single class of agents.

The Markovian framework and mean-field approximation have already been introduced in transportation research. One particular example is the Markov routing game of ride-hailing vehicles, where the cruising and searching behaviors are characterized as an MDP to maximize the total revenue of individual vehicles over the operational horizon~\cite[e.g.,][]{calderone2017markov,li2022congestion,zhang2023ride}. 
The decentralized rebalancing problem in ride-hailing shares a similar framework, though the objective is to maximize the total number of matched trips or the net profit~\cite[e.g.,][]{jiao2021real,shou2020reward,jusup2024safe}. 
Another closely related line of research is the Markovian traffic equilibrium (MTE), where the transition between nodes is characterized by a prescribed link choice probability (e.g., Logit model)~\citep{akamatsu1996cyclic,baillon2008markovian,oyama2022markovian}. 
The solution approaches to these Markov population games can be broadly classified into two types. 
When the system dynamics are explicitly specified, a common approach computes the mean-field equilibrium through an MSA-based policy iteration. 
Given a mean policy of the population, a forward pass propagates the population distribution over states, while a backward pass solves an MDP of the representative agent by dynamic programming.  
The mean policy is then updated toward the best-response policy obtained in the backward pass using the method of successive averages (MSA). 

This method, however, suffers from the scalability issue and often shows slow convergence~\cite[e.g.,][]{zhang2023ride,elokda2024dynamic,zhang2023routing}. On the other hand, MTE models often feature more computationally efficient solution algorithms thanks to their particular mathematical properties, which, however, do not hold for general MPG.
Moreover, when underlying dynamics of the Markov game are unknown, multi-agent reinforcement learning (MARL) and mean-field reinforcement learning (MFRL) are usually exploited to directly learn system dynamics, value function, or optimal policy from simulations.

In this paper, we propose a special case of Markov population games, namely, \textit{Markov population potential game (MPPG)}.
We focus on a non-cooperative game among homogeneous self-interested agents in a finite-horizon setting, where each agent's decisions are modeled as a finite-horizon MDP, and \emph{both the immediate reward and the state transition depend on the joint state-action distribution of all agents}---a particular expression of mean field.
Accordingly, the resulting Nash equilibrium is specified as a mean-field equilibrium. 
In particular, we establish sufficient conditions on stage rewards and transition dynamics under which a Markov population game possesses a potential function, thus becomes an MPPG. We further elaborate the construction of the potential function, such that the KKT points of the corresponding optimization problem coincide with the mean-field equilibria. 

Although the MPPG condition on stage reward is often satisfied in transportation applications, the condition on the transition dynamics does not hold in general. 
Hence, we proceed to develop theoretical results on equilibrium characterization when only the reward condition is satisfied. Specifically, we show that the same optimization problem can still be constructed as per MPPG, though its KKT points only yield $\varepsilon$-equilibria. We then derive conditions under which $\epsilon$ reduces to zero, and the priori and posteriori bounds of $\varepsilon$.

The optimization problem associated with MPPG and $\varepsilon$-MPPG also motivates a projected-gradient policy iteration algorithm (PG) that is shown to outperform MSA-based policy iteration in the numerical experiments with more stable convergence and smaller equilibrium gaps. The numerical experiments also demonstrate the difference between MPPG and $\varepsilon$-MPPG and compare the priori and posteriori bounds.

The remainder of the paper is organized as follows. Section \ref{sec:review} reviews related literature and positions our contributions. Section \ref{sec:prelim} introduces preliminaries on Markov games and Markov population games. 
Section \ref{sec:MPPG} formally introduces the MPPG framework, establishes sufficient conditions for the existence of a potential function, outlines a constructive procedure for potential functions, defines the $\varepsilon$-equilibrium, then proposes bounds on $\varepsilon$, and presents the solution algorithm. Section \ref{sec:connect_MTE} discusses connections to existing Markovian traffic equilibrium models.
Section \ref{sec:app} details the design of numerical experiments and reports the main findings. Finally, Section \ref{sec:conclusion} concludes this study and discusses potential future research directions.

\section{Related work}\label{sec:review}
\subsection{Markovian models in transportation}

Markovian models have been widely used in transportation research. A prominent line of applications is Markovian traffic
assignment. Its origin can be traced to \cite{dial1971probabilistic}, which represented path choice as a sequence of link
choices but restricted the choice set without cycles. Following studies developed Markovian
traffic-assignment formulations for more general networks~\citep[e.g.,][]{bell1995alternatives,akamatsu1996cyclic,akamatsu1997decomposition}. These fundamental works were subsequently generalized to include more general discrete choice models~\citep{baillon2008markovian}. In all these modeling frameworks, the state is defined as the current node, the action refers to the next link chosen, while the state transition is exogenous and deterministic (i.e., moving from the current node to the ending node of the selected link).
Besides, these models introduce a random error term to the link utility (immediate reward) following a prescribed probability distribution (e.g., i.i.d Type-I Extreme Value distribution with
zero mean) so that 
the link choice probabilities can be derived in a closed-form expression (e.g., the multinomial logit model), as well as the associated value function \citep{fosgerau2013link, oyama2022markovian}.
Recently, the perturbed utility Markovian choice model (PUMCM) has been proposed as a new class of Markovian route choice models, which generalizes perturbed utility theory to sequential decision problems, which is then embedded in an equilibrium model~\citep{yao2026perturbed}. 

Nevertheless, all studies in this line of research rely on exogenous and deterministic transitions, which limit their applicability to more general transportation problems.

Another stream of research is motivated by the competition among strategic agents in dynamic transportation systems. A representative example is the MDP routing game developed by~\cite{calderone2017infinite}, where each agent makes sequential routing decisions under an MDP, and the immediate reward depends on aggregate routing behavior. Subsequent studies have extended this framework to incorporate system-level constraints and objectives~\citep{li2019tolling}, as well as agent heterogeneity~\citep{li2022congestion}. 
In parallel, \cite{zhang2023ride} developed a game-theoretic framework to model the strategic routing behavior of ride-hailing vehicles. Same as in \cite{calderone2017infinite}, each vehicle's routing decisions follow an MDP. However, the interactions among vehicles are captured through the physical passenger-vehicle matching process, which in turn determines both the expected reward and state transitions. A follow-up study extended this framework to jointly model routing and charging decisions of electric ride-hailing vehicles, revealing congestion in both the matching and charging processes~\citep{zhang2023routing}.

The third group of research consists of MARL-based studies, which directly learn policies from simulation to tackle complex multi-agent decision-making problems. Applications of MARL in transportation span a wide range, from the operational problems in ride-hailing systems, such as trip dispatching, pricing, and empty vehicle rebalancing~\citep{lin2018efficient,chu2019multi,shou2022multi}, to traffic control with different control instruments, including traffic signals and autonomous vehicles~\citep{chen2020toward,vinitsky2023optimizing}. Although these learning-based approaches demonstrate promising performances in specific problems, they do not generalize into a unified modeling framework. Moreover, the computationally intensive training procedures substantially limit their integration into high-level system design problems.

\subsection{Markov potential games}
Markov potential games (MPGs) generalize static potential games to the Markovian setting~\citep{macua2018learning,leonardos2021global}. 
As in other potential games, a salient feature of MPGs is the existence of a potential function that turns equilibrium seeking into solving an optimization problem. Specifically, in an MPG, the change in any agent's value due to a unilateral policy shift equals the change in the potential function.  
Prior work on MPGs has primarily focused on developing efficient algorithms to compute equilibrium and establishing their convergence properties, both for discounted infinite horizon~\citep[e.g.,][]{fox2022independent,zhang2024gradient,maheshwari2025independent} and for finite horizon~\citep{mao2021decentralized, song2021can}.
For the existence and construction of potential functions, strong assumptions are often imposed. Examples include restricting the MPG to a cooperative game so that the total expected reward acts as the potential, or assuming independent state transitions so that the stage potentials sum up to be the global potential~\citep{macua2018learning,leonardos2021global}. These conditions, however, are over-restrictive for general transportation systems.

Since potential functions typically do not exist for general Markov games, recent research has shifted toward developing near-potential formulations. One notable example is $\alpha$-MPGs~\citep{guo2025markov,Guo2025Alpha}, where an approximate potential function is defined such that the absolute difference between the change in potential and the change in any agent's value is bounded by a constant $\alpha$. Under this formulation, an optimizer of the $\alpha$-potential function corresponds to an $\epsilon$-Nash equilibrium. 
It has further been shown that any Markov game with finite state and action spaces is an $\alpha$-MPG~\citep{guo2025markov}, although the approximation gap $\alpha$ could be arbitrarily large. 
Similar to that on MPGs, research on $\alpha$-MPGs has largely focused on establishing analytical properties and developing solution algorithms. In contrast, relatively little effort has been devoted to constructing the potential function or linking these formulations with real-world applications. 
The most recent extension to general dynamic games with a finite number of players sheds light on the conditions under which a potential function exists, as well as its general structure~\citep{Guo2025Alpha}, though the derived potential does not have a closed form to be easily utilized in applications. 

\indent The connection between equilibrium computation and optimization over potential has a long tradition in traffic assignment problems. 
With the standard separable link cost assumption, \citet{beckmann1956studies} developed the well-known Beckmann formulation. 
Similar convex optimization problems were subsequently developed for stochastic user equilibrium (SUE), where routing costs are subject to random errors~\citep{fisk1980some,daganzo1982unconstrained}. 
Extending from \citet{dial1971probabilistic}'s logit network-loading over a restricted path set, subsequent formulations accommodate cyclic routes without explicit path enumerations~\citep{bell1995alternatives, nesterov2007characteristic}.
The perturbed utility traffic equilibrium models further replace the random error with a convex perturbation of link flows that generalizes additive random utility models and allows more flexible substitution patterns and corner solutions~\citep{fosgerau2022perturbed, yao2024perturbed}.
Several Markovian traffic equilibrium (MTE) models also admit equivalent optimization formulations~\citep[e.g.,][]{akamatsu1997decomposition,baillon2008markovian,oyama2022markovian}, and their connections to our framework will be discussed in Section~\ref{sec:connect_MTE}.

\subsection{Population game and mean field game}

Population games and mean field games are two widely used modeling frameworks to address the complexity introduced by large numbers of players. 
Population games, proposed by \cite{sandholm2001potential,sandholm2010population}, consist of one or a finite number of homogeneous agent groups, each of which shares the same payoff function depending on their own strategy and the distribution of other groups' strategies.
The large population is treated as a continuum such that each agent's decision has a negligible impact on the system. Population games with a single homogeneous agent group are closely related to mean field games~\citep{huang2006large,lasry2007mean}, where each agent makes decisions based on their own state and on the aggregate mean-field rather than on direct pairwise interactions. 
Since both population games and mean field games represent collective agent behaviors, their equilibria are expressed by aggregate strategy profiles and mean policy, respectively, along with the resulting state-action flows and mean field dynamics~\citep{adlakha2015equilibria,neumann2020stationary}.  

The introduction of population games and mean field games greatly reduces the complexity of computing equilibrium in complicated large-scale multi-agent systems like transportation systems. For example, most classic traffic assignment models can be cast into the population-game framework~\citep[e.g.,][]{beckmann1956studies,sheffi1985urban}. 
More recently, \cite{cabannes2021solving} and \cite{tanaka2020linearly} formulate routing games and approximate congestion effects through a mean-field state distribution. 
\cite{huang2021dynamic} model the routing game among autonomous vehicles using mean field approach so that the complex multi-vehicle differential game is reduced to solving the standard Hamilton–Jacobi–Bellman and flow evolution equations. 
More closely related to this study, \cite{elokda2024dynamic} introduced the dynamic population game, or Markov population game, where each agent follows an MDP and their state dynamics depend on the joint policy profile. It is shown that each dynamic population game has an equivalent static population game, and its stationary equilibrium corresponds to the Nash equilibrium of its static counterpart. This framework has been applied to tackle fair resource allocation in traffic congestion management~\citep{elokda2025carma}.

\indent Substantial research efforts have been made recently to develop solution methods for mean-field equilibria. For example, fictitious play \citep{perrin2020fictitious} and online mirror descent
\citep{perolat2021scaling} have been developed with convergence guarantees when the rewards satisfy monotonicity conditions and the transition kernel is independent of population distribution. Yet, such conditions may not hold in many transportation applications. Besides, \cite{guo2023general} proposed both value- and policy-based reinforcement learning algorithms for general mean-field games, where both transition and reward depend on aggregate states and actions. Their convergence guarantees, however, rely on a contraction assumption on the composite best-response and population-update operator.   
Another approach connects equilibrium computation to 
optimization. For example, \citet{guo2024mf} propose an optimization formulation for general discrete-time mean-field games, whose objective measures the total violation of individual optimality and population consistency. 
Accordingly, solutions with zero objective value correspond to equilibria. A projected gradient descent algorithm is further developed that guarantees convergence to equilibria under suitable initializations, without monotonicity or contraction assumptions.
On the other hand, potential mean-field games naturally feature equivalent optimization problems. 
For mean-field games with discrete states and continuous time,
\citet{gueant2011infinity} constructed the global potential function,
thereby reformulating the coupled forward-backward equilibrium system
as a global optimal control problem. Subsequent studies explored Newton's iteration method \citep{gao2022master}
and a primal--dual Newton method \citep{cui2026convergence}
for computing equilibria. 
Differently, \citet{bonnans2023discrete} considered a potential mean-field game with discrete time, finite states, and continuous actions, and proposed primal-dual proximal methods and ADMM to solve it.
These models, however, assume the individual transitions are independent of the aggregate actions, different from the setting considered in this paper.

MPPG proposed in this paper lies at the intersection of Markov potential games, population potential games, and Markov population games. On the one hand, we extend the classic population potential games~\citep{sandholm2001potential} by incorporating Markovian dynamics over a finite horizon. To better reflect the physical dynamics in transportation systems, we allow state transitions to depend on the joint state and action distribution, which are often simplified in the literature~\citep[e.g.,][]{calderone2017infinite, li2019tolling, li2022congestion,leonardos2021global}. By lifting the Markov potential game to the population level, we are able to derive more tractable conditions to check whether a Markov game is also a potential game, and if so, how to construct its potential function. The closed-form potential function further enables us to develop a more efficient solution algorithm to compute the equilibrium than existing studies on Markov population games~\citep{elokda2024dynamic,elokda2025carma}.

\section{Preliminaries}\label{sec:prelim}
\subsection{Markov game}

A Markov game is defined as $\setG = (\setN, \setS, \setA, r, P, \gamma)$, where $\setN$ is the set of homogeneous agents with $|\setN| = N$;  
$\setS$  and $\setA$ are finite sets of individual states and actions, respectively, shared by all agents; $r: \setS^N \times \setA^N \times \setS \to \bbR$ is the immediate reward function depending on the states and actions of all agents; $P: \setS^N \times \setA^N \times \setS \to [0, 1]$ is the state transition function that also depends on the joint state-action profile of all agents; and $\gamma \in (0, 1]$ is the discount factor.
For each agent $n$, let $s_n \in \setS$ and $a_n \in \setA$ denote its state and action, and $(\vecs, \veca)$ represent the joint state–action profile of all agents. We further introduce $\vecs_{-n}$ and $\veca_{-n}$ to express the states and actions of all agents other than $n$. 
Then $r(s_n' \mid s_n, a_n, \vecs_{-n}, \veca_{-n})$ and $P(s'_n \mid s_n, a_n, \vecs_{-n}, \veca_{-n})$ give the immediate reward and transition probability, respectively, of agent $n$ associated with a state transition from $s_n$ to $s_n'$ under its own action $a_n$ and the joint state-action of other agents $(\vecs_{-n}, \veca_{-n})$.

We introduce an index $t$ to distinguish states and actions across different stages, denoted as $s_{n,t}$ and $a_{n,t}$. 
Each agent $n \in \setN$ makes sequential decisions following a state-dependent policy 
$\pi_n = \textrm{col}\{\pi_{n,t}\}_{t\in\setT}$ over a finite horizon $\setT = \{0,\dots, T-1\}$, 
where $\pi_{n,t}(s,a)$ denotes the probability of taking action $a$ at state $s$ and stage $t$. 
Given the initial state distribution $\mu_0$ and the joint policy of other agents $\vecpi_{-n}$, the expected total discounted reward of agent $n$ is
\begin{align}
V_n(\pi_n, \vecpi_{-n}) 
= \bbE\!\left[
\sum_{t=0}^{T-1} 
\gamma^t \, r(s_{n,t+1} \mid s_{n,t}, a_{n,t}, \vecs_{-n,t}, \veca_{-n,t})
\right],
\end{align}
where $s_{n,0} \sim \mu_0$, $a_{n,t} \sim \pi_{n,t}(s_{n,t}, \cdot)$, 
and $s_{n,t+1} \sim P(\cdot \mid s_{n,t}, a_{n,t}, \vecs_{-n,t}, \veca_{-n,t})$ 
for all $n \in \setN$ and $t \in \setT$.

\subsection{Markov population game}

When the set of agents $\mathcal{N}$ is sufficiently large, we may present $\setG$ as a Markov population game by defining the state distribution $\mu=\textrm{col}\{\mu_t\}_{t\in\{0,\dots, T\}}$ and mean policy $\pi=\textrm{col}\{\pi_t\}_{t\in\setT}$ as follows:
\begin{subequations}\label{eq:def_mu_pi}
\begin{align}
    \mu_t(s) &= \frac{1}{N}\sum_{n\in\setN} \bm{1}(s_{n,t} = s),\\
    \pi_t(s,a) &= \frac{1}{N}\sum_{n\in\setN} \pi_{n,t}(s,a). 
\end{align}
\end{subequations}

Let $\Delta(S)$ and $\Delta(A)$ denote probability distribution over
$\setS$ and $\setA$, respectively. 
Then, the spaces of state distributions and policies are defined as
\begin{align}
   \mathcal M = \prod_{t=0}^{T}\Delta(\setS), \; 
    \Pi = \prod_{t=0}^{T-1}\prod_{s\in S}\Delta(\setA).
\end{align}

We introduce the common assumption in population games and mean field games~\citep{sandholm2009large,lasry2007mean} that the influence of any agent on the system dynamics is negligible while its action, immediate reward, and state transition are significantly influenced by the aggregate state-action flows of all agents. Hence, we replace the joint states and actions $(\vecs,\veca)$ with the state distribution and policy $\mu,\pi$. Since all agents are homogeneous, in the following analysis, we consider a representative agent who adopts the mean policy $\pi$ subject to the state distribution $\mu$. Accordingly, at each stage $t$, the transition function and immediate reward function of the representative agent depend on $\mu,\pi$, rewritten as $P(s'|s, a, \mu_t, \pi_t)$ and $r(s'|s, a, \mu_t, \pi_t)$, respectively. We further define the stage-wise expected reward as
\begin{align}\label{eq:stage reward}
    R(s, a, \mu_t, \pi_t) = \sum_{s'\in \setS} r(s'|s,a,\mu_t, \pi_t) P(s'|s, a, \mu_t, \pi_t), \quad \forall s\in\setS,a \in\setA.
\end{align}
Besides, the dynamics of the state distribution are given by 
\begin{align}\label{eq:state_dynamics}
    \mu_{t+1}(s') = \sum_{s\in \setS}\sum_{a\in\setA} P(s'|s, a, \mu_t, \pi_t) \mu_t(s) \pi_t(s,a). 
\end{align}
With the expected reward Eq.~\eqref{eq:stage reward}, the value function from time $t$ and state $s$ is defined as 
\begin{align}
     V_t(s,\mu,\pi) = \bbE\left[\sum_{\tau=t}^{T-1} \gamma^{\tau-t} R(s,a,\mu_{\tau},\pi_{\tau})\right]
\end{align}
where $a\sim\pi_{t}(s,\cdot)$. Equivalently, it can be expressed as
\begin{align}\label{eq:def_V}
    V_t(s,\mu,\pi) = \sum_{a\in\setA} \pi_t(s,a)Q_t(s,a,\mu,\pi)
\end{align}
where the state-action value function $Q_t(s,a,\mu,\pi)$ satisfies the recursive relation:
\begin{align}\label{eq:q_value_def}
    Q_t(s,a,\mu,\pi) = R(s,a,\mu_t,\pi_t)+ \gamma \sum_{s'\in\setS}P(s'|s,a,\mu_t,\pi_t)V_{t+1}(s',\mu,\pi).
\end{align}

We then define the equilibrium of Markov population games.
\begin{definition}\label{def:equilibrium_BR}
    The equilibrium of a Markov population game is defined as a tuple $(\mu^*,\pi^*)$, such that the following conditions hold for any $ t\in \setT$: 
    \begin{subequations}\label{eq:def_equ_BR}
        \begin{align}
            &\pi^*_t(s,\cdot) \in\arg\max_{\pi_t(s,\cdot)}\sum_{a\in\setA}\pi_t(s,a)Q_t(s,a,\mu^*,\pi^*)\quad &\forall s\in\setS, \label{eq:def_equ_pi_BR}\\
            &\mu^*_{t+1}(s') = \sum_{s\in\setS}\sum_{a\in\setA} P(s'|s,a,\mu^*_t,\pi^*_t)\mu^*_t(s)\pi^*_t(s,a), \quad &\forall s'\in \setS.\label{eq:def_equ_mu_BR}
        \end{align}
    \end{subequations}
\end{definition}
\noindent Condition~\eqref{eq:def_equ_pi_BR} states that given the mean field $(\mu^*,\pi^*)$, a representative agent selects a policy to maximize its value function. At equilibrium, the mean policy $\pi^*$ coincides with such a best response at every state and time. Condition~\eqref{eq:def_equ_mu_BR} then ensures the evolution is consistent with the mean-field dynamics. 
It is worth noting that, in the finite-horizon setting considered in this paper, the mean-field state distribution can be uniquely induced by a given mean policy and an initial distribution $\mu_0$, though we express them separately for the clarity of notation.
The following corollary gives an equivalent expression of the equilibrium.
\begin{corollary}
\label{def:equilibrium}
    A tuple $(\mu^*,\pi^*)$ is an equilibrium of a Markov population game if for any $ t\in \setT$, it satisfies the following conditions: 
    \begin{subequations}\label{eq:def_equ}
        \begin{align}
            &\pi_t^*(s,a)\Big(Q_t(s,a,\mu^*,\pi^*)-Q^*_{t}(s,\mu^*,\pi^*)\Big) = 0,\quad &\forall s\in\setS, a\in\setA,\label{eq:def_equ_Q}\\
            &\mu^*_{t+1}(s') = \sum_{s\in\setS}\sum_{a\in\setA} P(s'|s,a,\mu^*_t,\pi^*_t)\mu^*_t(s)\pi^*_t(s,a), \quad &\forall s'\in \setS,\label{eq:def_equ_mu}
        \end{align}
    \end{subequations}
    where $Q^*_t(s,\mu,\pi) = \max_{a\in\setA} Q_t(s,a,\mu,\pi)$ is the optimal Q-value at time $t$ and state $s$.  
\end{corollary}
\noindent Condition~\eqref{eq:def_equ_mu} is consistent with Condition~\eqref{eq:def_equ_mu_BR}, while Condition~\eqref{eq:def_equ_Q} is derived from the optimality condition of \eqref{eq:def_equ_pi_BR}. Specifically, the complementarity condition states that either $\pi_t^*(s,a)=0$ or $Q_t(s,a,\mu^*,\pi^*)$ yields the optimal state-action value.

The existence of an equilibrium relies on the continuity of the reward function and transition kernel with respect to the mean-field variables, which is formally established below. 
\begin{theorem}\label{prop:existence_equ}
Suppose that for all $t\in\mathcal T$, $s,s'\in\mathcal S$, and $a\in\mathcal A$,
the immediate reward and state transition kernel, and hence the expected stage reward,
are continuous with respect to $(\mu_t,\pi_t)$.
Then the Markov population game defined above admits at least one equilibrium
$(\mu^*,\pi^*)$ satisfying Eqs.~\eqref{eq:def_equ_Q} and ~\eqref{eq:def_equ_mu}.
\end{theorem}
\begin{proof}
    The proof largely follows \cite{elokda2024dynamic} and is detailed in Appendix \ref{app:existence_equ}.
\end{proof}

\section{Markov population potential game}\label{sec:MPPG}

We now connect the Markov population game to the Markov potential game (MPG). 
Here, we focus on a class of Markov population games with rewards and transition probabilities both defined on
the population mass of each state-action pair. 
Specifically, we define the \textit{state-action flow} as $x_t(s,a) = \mu_t(s )\pi_t(s, a)$ and collectively $x(\mu,\pi) = \{x_t(\mu,\pi)\}_{t\in\setT}$. 
Then, given the initial state distribution $\mu_0$, the feasible set of state-action flows is given by
\begin{align}
\Omega=
\Bigg\{
&x \Bigg|\
x_t(s,a)\geq0, \sum_{a\in\mathcal A}x_0(s,a)=\mu_0(s), \nonumber\\ &\sum_{a\in\mathcal A}x_{t+1}(s,a)
=
\sum_{s'\in\mathcal S}\sum_{a'\in\mathcal A}
P(s\mid s',a',x_t)x_t(s',a')
&\forall s\in\setS,\ a\in\setA,\ t\in\setT\setminus\{T-1\}
\Bigg\}.
\label{eq:feasible_flow_set}
\end{align}
Accordingly, the mean state distribution and policy can be recovered from a feasible flow $x\in\Omega$:
\begin{align}
\mu_t(s)&=\sum_{a\in\mathcal A}x_t(s,a),&& \forall s\in\setS, t\in\setT,\\ 
\pi_t(s,a)&=\frac{x_t(s,a)}{\mu_t(s)},&& \forall a\in\setA, s\in\{s'\in\setS|\mu_t(s') >0\}, t\in\setT. 
\end{align}
When $\mu_t(s)=0$, the policy $\pi_t(s,\cdot)$ is not uniquely determined by the state-action flow and can be specified separately without affecting the induced flow.

We further introduce the following assumptions about the reward and state transition. 
\begin{assumption}\label{ass:regularity_diff}
The immediate reward $r$ and state transition kernel $P$ are functions of state-action flows $x$ and satisfy the following properties:
For any $t \in \mathcal{T}$, $s,s' \in \mathcal{S}$, $a\in\setA$,
\begin{itemize}
    \item $r(s'\mid s,a,x_t)$ and $P(s' \mid s,a,x_t)$, consequently $R(s,a,x_t)$, are continuously differentiable with respect to $x_t$;

    \item $\nabla P(s'\mid s,a,x_t)$ is continuous in $x_t$ and
    uniformly bounded. Specifically, there exists a constant
    $L_P<\infty$ such that
    \begin{equation}\label{eq:tran_grad_bound}
    \left\|
    \nabla P(s'\mid s,a,x_t)
    \right\|_\infty
    \leq L_P.
    \end{equation}
\end{itemize}
\end{assumption}
Since reward and state transition are both functions of state-action flow $x$, the value functions $V$ and $Q$ can also be written as functions of $x$. The corresponding notations will be used hereafter without being explicitly defined.

With Assumption~\ref{ass:regularity_diff}, we now define the Markov population potential game as follows:

\begin{definition}\label{def:MPPG}
    A Markov population game is a Markov population potential game (MPPG) if there exists a continuously differentiable function $\Phi(x)$ such that
    \begin{align}\label{eq:MPPG}
        \frac{\partial\Phi(x)}{\partial x_t(s,a)} = \gamma^t Q_t(s, a,x),\; \forall t\in \setT,s\in\setS,a\in\setA
    \end{align}
    Here, differentiation with respect to $x_t(s,a)$ accounts for the induced changes in subsequent state-action flows through the population dynamics, while holding the continuation policies fixed.
\end{definition}

\noindent In the definition of a finite-player MPG~\citep{macua2018learning,leonardos2021global}, the potential function is defined such that the change in an individual agent's value function due to a unilateral policy shift equals the change in the potential function. However, in the MPPG framework, as the population size approaches infinity, the impact of a single agent on system dynamics becomes negligible. Instead, we focus on the representative agent whose expected total discounted reward depends on the state-action flows. Consequently, the potential function must capture the marginal impact of these flows on the reward. Specifically, it requires that the derivative of the potential with respect to the flow at each state-action pair $(s, a)$ equals the corresponding discounted Q-value. 

\begin{remark}\label{remark:equ_exist}
Assumption~\ref{ass:regularity_diff} implies the continuity of immediate reward and state
transition kernel with respect to $(\mu_t,\pi_t)$ for all $t \in \mathcal{T}$, $s,s' \in \mathcal{S}$, and $a \in \mathcal{A}$. Therefore, according to Theorem~\ref{prop:existence_equ}, there exists an equilibrium defined in Definition~\ref{def:equilibrium}.
\end{remark}

\subsection{Sufficient conditions for the existence of potential}

We next present sufficient conditions for guaranteeing the existence of a potential function, along with a constructive procedure for the potential. 

\begin{assumption}\label{ass:sym_reward}
    The expected stage reward satisfies the symmetric condition:
\begin{align}\label{eq: sym_R}
    \frac{\partial R(s,a,x_t)}{\partial x_t(s',a')} = \frac{\partial R(s',a',x_t)}{\partial x_t(s,a)}, \quad \forall t\in \setT,s,s'\in\setS, a,a'\in\setA.
\end{align}
\end{assumption}
\noindent Assumption \ref{ass:sym_reward} guarantees a stage-wise potential function, like in static population potential games~\citep{sandholm2001potential, sandholm2009large}. Specifically, it states that, when cross-strategy effects in a static population game are symmetric, there is a potential function and each strategy’s payoff corresponds to the marginal change in the potential with respect to the aggregate flow of such strategy. 
However, a Markov game is not necessarily a potential game even if each stage is a potential game~\citep{leonardos2021global,zhang2024gradient}. 
In this study, we derive the following additional condition for a Markov population game to possess a potential.

\begin{assumption}\label{ass:transition}
The expected marginal impact of flow $x_t(s,a)$ on the transitions to state $s'$ is zero, i.e., $\forall t\in\setT, s,s'\in\setS, a\in\setA$, 
\begin{align}\label{eq:expect_trans_marg}
    \mathbb{E}\left[\frac{\partial P(s'|\cdot, x_t)}{\partial x_t(s,a)}\right] = \sum_{s''\in\setS}\sum_{a''\in\setA} \frac{\partial P(s'|s'', a'', x_t)}{\partial x_t(s,a)}x_t(s'',a'') = 0.
\end{align}
\end{assumption}
\noindent Assumption~\ref{ass:transition} requires that the expected marginal impact of any state-action flow on the transition dynamics vanishes, regardless of the flow-dependent transition. 
In other words, this condition rules out systematic first-order effects of flows on state transitions, while allowing for nontrivial dependence of transition kernels on the mean field. 

When Assumptions \ref{ass:sym_reward} and \ref{ass:transition} both hold, we prove a global potential function exists and can be constructed as the discounted sum of the stage-wise potential functions.
This result is formally stated in the following theorem:
\begin{theorem}\label{prop:suffic_potential}
A Markov population game is an MPPG if it satisfies Assumptions~\ref{ass:regularity_diff}--\ref{ass:transition}.
The corresponding potential function is constructed as 
\begin{align}\label{eq:potential_construct}
\Phi(x)=\sum_{t=0}^{T-1} \gamma^t\phi(x_t), 
\end{align}
where $\phi$ is the stage potential such that for any $t\in \setT, s\in\setS, a\in\setA$, 
\begin{align}\label{eq:stage_potential}
    \frac{\partial \phi(x_t)}{\partial x_t(s,a)} = R(s,a,x_t).
\end{align}
\end{theorem}

\begin{proof}
    We prove by first showing Condition \eqref{eq:MPPG} holds at the terminal stage and then completing by backward induction. More details are provided in Appendix~\ref{app:prop_suffic}.
\end{proof}

Next, we proceed to establish two corollaries of sufficient conditions that are induced from Assumption~\ref{ass:transition}, and the proof is shown in Appendix~\ref{app:coro_transition}. 
\begin{corollary}\label{coro:fully-independent}
    A Markov population game is an MPPG if it satisfies Assumptions~\ref{ass:regularity_diff}--\ref{ass:sym_reward}
    and the state transition function $P(s'|s,a,x_t)$ is essentially independent of $x_t$, that is, 
    \begin{align}\label{eq:indep-transition}
        \frac{\partial P(s'|s'',a'',x_t)}{\partial x_t(s,a)} = 0, \; \forall t\in\setT, s,s',s''\in\setS, a,a''\in\setA.
    \end{align}
\end{corollary}

\begin{corollary}\label{coro:symm-independent}
    A Markov population game is an MPPG if 
    it satisfies Assumptions~\ref{ass:regularity_diff}--\ref{ass:sym_reward}
    and the state transition Jacobian $\nabla P(s'|s,a,x_t)$ is symmetric with the directional derivative satisfying
    \begin{align}\label{eq:symm-independent}
        \nabla P(s'|s,a,x_t)^\top x_t = 0, \; \forall t\in\setT,s,s'\in\setS,a\in\setA.
    \end{align}
\end{corollary}

Note that the condition in Corollary~\ref{coro:fully-independent} implies that state transitions are fully independent among state-action pairs, in line with classic Markovian traffic equilibrium models~\citep{baillon2008markovian,oyama2022markovian}.
However, it imposes strong constraints on state transitions compared to Assumption~\ref{ass:transition}. 
The condition in Corollary~\ref{coro:symm-independent} is less restrictive but also imposes additional symmetric structure on the marginal effects of state-action flows on transitions. Nevertheless, it is more interpretable and easier to verify, i.e., the aggregate marginal impacts vanish for any state transition. 

\begin{remark}\label{remark:other_suffic_MPG}
Other sufficient conditions for potential existence have been proposed, mostly for finite-player games. \cite{leonardos2021global} assume action-independent state evolution, similar to Corollary~\ref{coro:fully-independent}, which is restrictive in most applications. \cite{mguni2021learning} introduce state transitivity, requiring that for any fixed joint action, the reward difference from changing state is identical for all agents, which is also restrictive and difficult to verify in practice. \cite{zhang2024gradient} impose three independence conditions: (i) the global state can be decomposed into agent-level components; (ii) each agent’s policy depends only on its own local state; and (iii) individual state transitions depend only on local states and actions. \cite{brock2025coordinating} allow shared state components but require that changes in an agent's local state and action leave the marginal transition of the other agents' collective states unaffected. Assumption~\ref{ass:transition} adopts a similar principle at the population level:  any state-action flow has negligible influence on system state evolution.

For $\alpha$-MPGs, exact potential existence corresponds to the special case of $\alpha=0$~\citep{guo2025markov}. For more general dynamic $\alpha$-potential games, \cite{Guo2025Alpha} establish potential existence under the symmetry of the second-order linear derivatives of agents' value functions with respect to their own and other agents' actions. This condition is difficult to verify since the value function is recursively defined through complex rewards and state transitions. By contrast, our conditions on stage rewards (Assumption~\ref{ass:sym_reward}) and transitions (Assumption~\ref{ass:transition}) are easier to verify and apply.
\end{remark}

\subsection{Equivalent optimization problem}\label{sec:equ_optm}

This section presents the equivalent optimization problem for solving the equilibrium of MPPG, whose objective is the potential function as per Definition \ref{def:MPPG}. The result is formally stated in the following theorem.

\begin{theorem}\label{prop:optm_MPPG}
Every tuple $(\mu^*,\pi^*)$ satisfying the Karush–Kuhn–Tucker (KKT) conditions of the following optimization problem corresponds to an equilibrium of the MPPG with potential function $\Phi$ as per Eq.~\eqref{eq:potential_construct}:
\begin{subequations}\label{eq:optm_MPPG}
    \begin{align}
        \max_{\mu,\pi} \quad & \Phi(x(\mu,\pi)),\\
        s.t.\quad  &  \mu_{t+1}(s') = \sum_{s\in \setS}\sum_{a\in\setA} P(s'|s, a, x_t(\mu, \pi)) \mu_t(s)\pi_t(s,a),& \forall s'\in\setS, t\in\setT, \label{eq:optm_MPPG_dynm}\\ 
        & \sum_{a\in \setA} \pi_t(s,a)=1, & \forall s\in\setS, t\in\setT, \label{eq:optm_MPPG_simplex}\\
        & \pi_t(s,a)\geq 0, & \forall s\in\setS, a\in\setA, t\in\setT.\label{eq:optm_MPPG_nonneg}
    \end{align}
\end{subequations}
\end{theorem}

\begin{proof}
     We show that the KKT conditions correspond to the equilibrium conditions \eqref{eq:def_equ}. Specifically, at a KKT point, the dual variables associated with Constraints~\eqref{eq:optm_MPPG_dynm} and~\eqref{eq:optm_MPPG_simplex} correspond to the state value and optimal Q-value, respectively. The detailed proof is provided in Appendix~\ref{app:prop_optm_MPPG}.   It is also easy to prove a similar result under the conditions in Corollaries~\ref{coro:fully-independent} and \ref{coro:symm-independent}.
\end{proof}

Here, the composite objective $\Phi(x(\mu, \pi))$ explicitly captures the dependence of the potential function on the flow, while the constraints are retained in terms of $(\mu, \pi)$ to preserve the original equilibrium structure.
And we note that, although in line with the classic result of potential games, the proof is not trivial due to the dependence of transition probabilities on the state–action flows. As shown in Appendix~\ref{app:prop_optm_MPPG}, this dependence introduces an extra expected marginal impact term in the KKT conditions, which greatly complicates the extension of static potential games to the more general dynamic setting.

\begin{remark}\label{remark:unique_equ}
   Generally, the objective of Problem~\eqref{eq:optm_MPPG} may be nonconcave, and flow-dependent transitions may yield nonlinear equality constraints, so equilibrium uniqueness is not guaranteed. A sufficient condition for uniqueness is that the transition kernel is exogenous and the expected stage reward is strictly monotone decreasing with respect to the state-action flow. The latter ensures strict concavity of each stage potential and their positively weighted sum $\Phi$, while exogenous transitions yield affine equality constraints and a convex feasible flow set. Thus, the problem admits a unique optimal, or equivalently, equilibrium flow. Such conditions are also widely employed in classical Markovian traffic assignment models\citep{baillon2008markovian,oyama2022markovian}. 
\end{remark}

\subsection{$\varepsilon$-equilibrium and approximation bounds}

We first introduce an approximate equilibrium for Markov population games motivated by the concept of restricted equilibrium \citep{mertikopoulos2018riemannian}. Specifically, we define $\varepsilon$ as the maximum difference in Q-values of two active actions, formally defined as follows. 
For convenience, we use the terms “support-restricted \(\varepsilon\)-equilibrium” and “\(\varepsilon\)-equilibrium” interchangeably hereafter.

\begin{definition}\label{def:epsilon_equ}
A tuple $(\tilde\mu,\tilde\pi)$ is a support-restricted $\varepsilon$-equilibrium of a Markov population game if there exists an $\varepsilon\ge 0$, such that:
\begin{subequations}
    \begin{align}
        &\vert{}Q_t(s,a_1,\tilde\mu,\tilde\pi) - Q_t(s,a_2,\tilde\mu,\tilde\pi)\vert{} \le \varepsilon,\quad \forall t \in \setT, s\in \setS:\tilde \mu_t(s)>0, a_1, a_2 \in \setA_+(t,s;\tilde\pi)\label{eq:eps_equ_BR}\\
        &\tilde\mu_{t+1}(s') = \sum_{s \in \setS} \sum_{a \in \setA} P(s'\vert{}s,a,\tilde\mu_t,\tilde\pi_t) \tilde\mu_t(s) \tilde\pi_t(s,a), \quad \forall s' \in \setS, t \in \setT\label{eq:eps_equ_mu}
        \end{align}
\end{subequations}
where $\setA_+(t,s;\tilde\pi) = \{a \in \setA\mid\tilde\pi_t(s,a) > 0\}$ denotes the set of active actions with positive policy.  
Eq.~\eqref{eq:eps_equ_mu} is consistent with Condition~\eqref{eq:def_equ_mu_BR}, and Eq.~\eqref{eq:eps_equ_BR} indicates that the reward improvement by switching between any two active actions is no larger than $\varepsilon$. 
\end{definition}

\begin{figure}[htb]
    \centering
    \begin{subfigure}
    {0.48\textwidth}
    \centering
    \includegraphics[width=\textwidth]{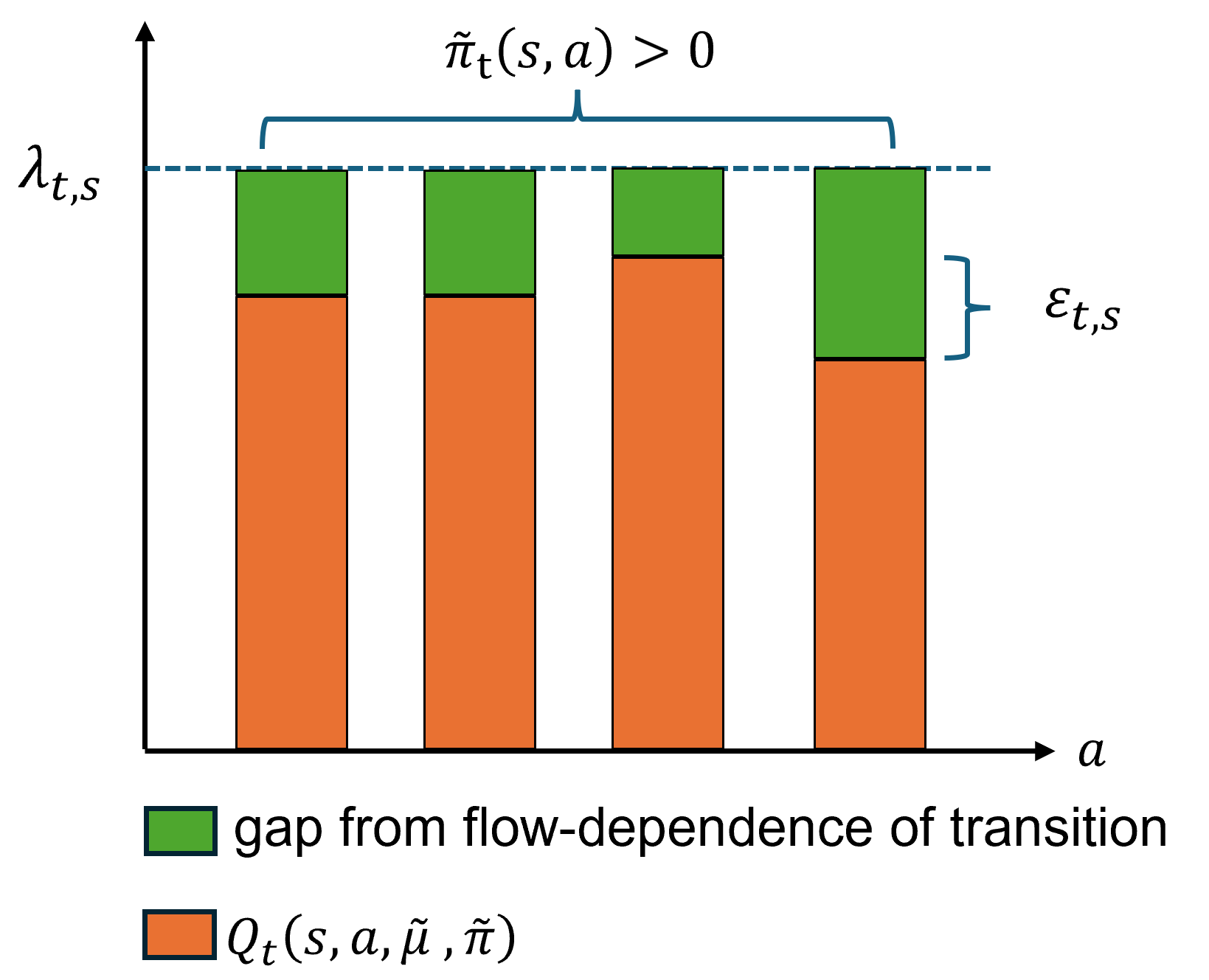} \caption{$Q$-values differ across active actions}
    \label{fig:epsilon_illu}
    \end{subfigure}
    \begin{subfigure}{0.48\textwidth}
    \centering
    \includegraphics[width=\textwidth]{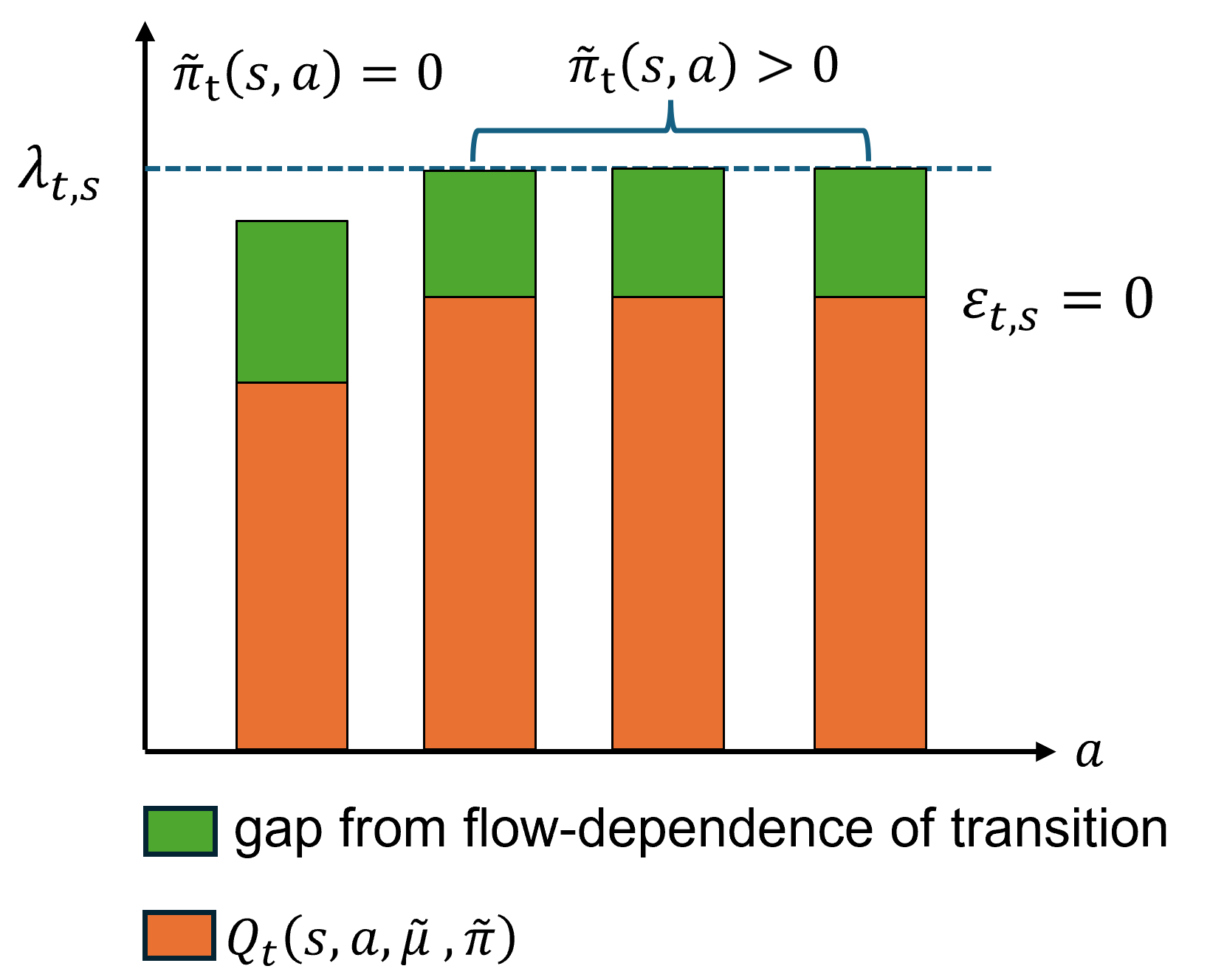} \caption{ Active actions share the maximum $Q$-value}
    \label{fig:fixgap_illu}
        
    \end{subfigure}
    \caption{Illustration of 
    $Q$-values  at KKT points}
    \label{fig:epsilon_fixgap}
\end{figure}

As shown in Appendix~\ref{app:prop_optm_MPPG}, at a KKT point $(\tilde\mu,\tilde\pi)$, the KKT conditions equalize the marginal action values across active actions.
However, when the transition kernel depends on the state-action flow, these marginal values need not coincide with the agent's Q-values.
In particular, the dual variable $\lambda_{t,s}$ need not coincide with the optimal Q-value $Q_t^*(s,\tilde\mu,\tilde\pi)$, and the Q-values of active actions may differ, as illustrated in Fig.~\ref{fig:epsilon_illu}. 
Even if the gap is constant among all active actions, the KKT point does not necessarily coincide with the equilibrium, unless an additional condition holds for the $\lambda_{t,s}$ and Q-values of inactive actions; see Fig.~\ref{fig:fixgap_illu} for an example.

Nevertheless, any KKT point can be shown to be an $\varepsilon$-equilibrium. A sufficient condition can also be derived under which a KKT point also corresponds to an equilibrium. These results are formally stated in the following proposition and proved in Appendix~\ref{app:equ_identical_Q}.

\begin{proposition}\label{prop:equ_identical_Q}
Suppose Assumptions~\ref{ass:regularity_diff} and
\ref{ass:sym_reward} hold. Then a KKT
point $(\tilde \mu, \tilde \pi)$ of Problem~\eqref{eq:optm_MPPG} corresponds to an $\varepsilon$-equilibrium of the Markov population game and further the exact equilibrium if for any $t\in\setT$, there exists a constant $\bar \epsilon_{t}$ such that  
\begin{align}\label{eq:local_stationarity}
\bar \epsilon_{t} = \sum_{s'\in\setS}\mathbb{E}\!\left[\frac{\partial P(s'\mid\cdot,\tilde x_t)}{\partial x_t(s,a)}\right]\nu_{t+1,s'},
&&\forall s\in\setS,\ a\in\setA,
\end{align}
where $\tilde x$ is the state-action flow induced by $(\tilde \mu, \tilde \pi)$ and $\nu_{t+1,s'}$ is the dual variable associated with Constraint~\eqref{eq:optm_MPPG_dynm}. 
\end{proposition}

Although Proposition~\ref{prop:equ_identical_Q} offers conditions to verify whether a KKT point of Problem~\eqref{eq:optm_MPPG} coincides with an equilibrium, it does not provide any insight into the approximation error $\epsilon$ when the conditions fail. 
Hence, in the following, we derive bounds on $\varepsilon$, both before and after solving Problem~\eqref{eq:optm_MPPG}. Specifically, the a priori bound tells whether it is valuable to approximate the equilibrium by constructing and solving the corresponding optimization problem, while the a posteriori bound helps efficiently evaluate the solution quality.

To derive the a priori bound, we first introduce the Dobrushin contraction coefficient of the transition
kernel \citep{gaubert2015dobrushin,de2023dual}, defined as the maximum
total variation in state transition probabilities across state-action pairs
\begin{align}
\theta=\max_x \max_t 
\max_{\substack{
s_1,s_2\in\setS\\
a_1,a_2\in\setA
}}
\frac{1}{2}\sum_{s'\in\setS}
\left|
P(s'\mid s_1,a_1,x_t)
-
P(s'\mid s_2,a_2,x_t)
\right|.
\label{eq:dobrushin_coefficient}
\end{align}
The maximum is well defined due to the finite horizon, states and actions, along with the compact set of $x$.  
By construction, $\theta\in[0,1]$, and a smaller value indicates less variation in state transition probabilities
across state-action pairs, thus stronger mixing of the transition kernel. In particular, $\theta=0$ when all state-action pairs have identical transition probabilities, while $\theta=1$ when two state-action pairs have no common next state with positive transition probability.

Besides, Assumption~\ref{ass:regularity_diff}  and the compactness of $\Omega$
ensure that the immediate reward $r$ is uniformly bounded, allowing us to define its uniform span as:
\begin{equation}
\label{eq:uniform_reward_span}
\Delta r =\max_x \max_t 
\left[
\max_{\substack{s,s'\in\setS\\a\in\setA}} r(s'\mid s,a,x_t)
-
\min_{\substack{s,s'\in\setS\\a\in\setA}} r(s'\mid s,a,x_t)
\right].
\end{equation}
With $\theta$ and $\Delta r$ defined above, we first derive a priori bound as stated in the following proposition, whose proof is provided in Appendix~\ref{app:theo_priori_bound}.

\begin{proposition}
\label{theo:priori bound}
Suppose Assumptions~\ref{ass:regularity_diff}
and~\ref{ass:sym_reward} hold, with $T\geq 2$.
Then every KKT point $(\tilde\mu,\tilde\pi)$
of Problem~\eqref{eq:optm_MPPG}
corresponds to an $\varepsilon$-equilibrium,
where $\varepsilon \leq \varepsilon_{\mathrm{pri}}$ and
\begin{equation}
\label{eq:epsilon_apriori_cases}
\varepsilon_{\mathrm{pri}}
= \Delta r \min\left\{
\sum_{k=0}^{T-1}(\gamma\theta)^k,\;
\sum_{k=1}^{T-1}\gamma^k
\left[
(\theta+|\setS|L_P)^k-\theta^k
\right]
\right\}.
\end{equation}
\end{proposition}
\noindent We note that, when the transition kernel is independent of the state-action flow, $L_P=0$. Consequently, the second term in the minimum vanishes, yielding $\varepsilon_{\mathrm{pri}}=0$.

The priori bound developed in Proposition~\ref{theo:priori bound}, however, can be overly conservative; see numerical results in Section~\ref{subsec:numerical:bounds}. Hence, we proceed to derive a posteriori bound constructed directly using the primal and dual solutions at a KKT point of Problem~\eqref{eq:optm_MPPG}.
Following Appendix~\ref{app:prop_optm_MPPG}, we denote 
$\nu_{t,s}$ and $\lambda_{t,s}$ as the
multipliers associated with constraints~\eqref{eq:optm_MPPG_dynm} and \eqref{eq:optm_MPPG_simplex} at time $t$ and state $s$, respectively. 
We then define the pseudo-Bellman residual constructed at a KKT point $(\tilde\mu, \tilde\pi)$, with induced flow $\tilde x$, for each state-action pair as 
\begin{equation}
\label{eq:posteriori_residual}
\mathcal{B}_t(s,a)
=\lambda_{t,s}-\left(R(s,a,\tilde x_t) + \gamma\sum_{s'\in\setS} P(s'\mid s,a,\tilde x_t)\nu_{t+1,s'}\right).
\end{equation}
It is worth noting that, when the KKT conditions coincide with the equilibrium, $\nu_{t,s}$ and $\lambda_{t,s}$ correspond to the optimal state values and optimal Q-values, thus $\mathcal{B}_t(s,a) = 0$ at actions with positive probability.
The terminal condition $\nu_T=0$ in Appendix~\ref{app:prop_optm_MPPG} also implies
$\mathcal{B}_{T-1}(s,a)=0$.

With the state-action residuals, we further define the following measures:
\begin{align}
\mathcal{B}_{\mathrm{act}}
&=
\max_{\substack{
t\in\setT\\
s:\,\tilde\mu_t(s)>0
}}
\left[
\max_{a\in\mathcal A_+(t,s;\tilde\pi)}\mathcal{B}_t(s,a)
-
\min_{a\in\mathcal A_+(t,s;\tilde\tilde\pi)}\mathcal{B}_t(s,a)
\right],
\label{eq:posteriori_C_act}
\\
\mathcal{B}_{\mathrm{avg}}
&=
\max_{t\in\setT}
\left[
\max_{\substack{s:\tilde\mu_t(s)>0}}\overline{\mathcal{B}}_t(s)
-
\min_{\substack{s:\tilde\mu_t(s)>0}}\overline{ \mathcal{B}}_t(s)
\right],
\label{eq:posteriori_C_avg}
\end{align}
where $\overline{ \mathcal{B}}_t(s)=
\sum_{a\in\mathcal A_+(t,s;\tilde\pi)}
\tilde\pi_t(s,a)\mathcal{B}_t(s,a)$
is the expected residual.
Specifically, $\mathcal{B}_{\mathrm{act}}$ measures the largest residual difference
among active actions at the same state with a positive share of agents, while
$\mathcal{B}_{\mathrm{avg}}$ measures the largest expected residual difference across states with positive shares.

The posteriori bound is derived as per the following proposition, with the proof detailed in Appendix~\ref{app:posteriori_bound}.

\begin{proposition}
\label{theo:posteriori_bound}
Under the same conditions as Proposition~\ref{theo:priori bound},
every KKT point $(\tilde\mu,\tilde\pi)$ of
Problem~\eqref{eq:optm_MPPG} corresponds to an
$\varepsilon$-equilibrium, where
$\varepsilon \leq \varepsilon_{\mathrm{post}}$ and
\begin{equation}
\label{eq:epsilon_aposteriori}
\varepsilon_{\mathrm{post}}
=
\min\left\{
\varepsilon_{\mathrm{pri}},\;
\mathcal{B}_{\mathrm{act}}
+
\mathcal{B}_{\mathrm{avg}}
\sum_{k=1}^{T-2}(\gamma\theta)^k
\right\}.
\end{equation}
\end{proposition}

\subsection{Solution algorithm}\label{sec:mppg_algo}

The classic approach to solving Markov population games is MSA-based policy iteration~\citep{elokda2024dynamic,zhang2023ride}, which could suffer from slow convergence in applications with complex dynamics~\citep[e.g.,][]{zhang2023routing}.
Taking advantage of the equivalent optimization problem in Section~\ref{sec:equ_optm}, we develop a projected-gradient policy iteration (PG) algorithm with a state-dependent adaptive step size, and show its performance in numerical experiments; see Section~\ref{sec:app} for more details. 

We first solve the following \textit{state evolution subproblem} given a policy $\pi$:
\begin{subequations}\label{eq:optm_sub}
\begin{align}
\max_{\mu}\; & \Phi(x(\mu,\pi)) \\
\text{s.t.} \quad 
& \mu_{t+1}(s') = \sum_{s,a} P(s'|s,a,x_t)\pi_t(s,a)\mu_t(s), \quad \forall s',t. 
\end{align}
\end{subequations}
The KKT conditions of the subproblem are given by 
\begin{align}
\nu_{t,s} &= \sum_{a\in\setA}\pi_t(s,a)
\Big[R(s,a,x_t) + \gamma \sum_{s'\in\setS}
\left(P(s'|s,a,x_t) + \mathbb{E}\!\left[\frac{\partial P(s'|\cdot,x_t)}{\partial x_t(s,a)}\right]\right)\nu_{t+1,s'}\Big], \label{eq:KKT_grad_mu_sub}\\
\mu_{t+1}(s') &= \sum_{s\in\setS} \sum_{a\in\setA}P(s'|s,a,x_t)\pi_t(s,a)\mu_t(s), \label{eq:KKT_eq_mu_sub}
\end{align}
which yields unique solution $\mu, \nu$; see Lemma~\ref{lemma:unique_V} in Appendix~\ref{app:Lemma_unique_V}. Here, $\nu$ is obtained through a linear system under a fixed mean field, which is closely related to the idea in \cite{guo2024mf}.

Let $\mu^\pi$ and $\nu^\pi$ denote the solution to subproblem \eqref{eq:optm_sub} at policy $\pi$. We then define the master problem as 
\begin{subequations}\label{eq:optm_mas}
    \begin{align}
        \max_\pi \quad & F(\pi) \\
        s.t.\quad  & \sum_{a\in\setA} \pi_t(s,a) = 1, &\forall s,t,\\
        & \pi_t(s,a) \geq0, & \forall s,a,t,
    \end{align}
\end{subequations}
where the objective function $F(\pi)$ is derived by relaxing only the constraint on state distribution dynamics:
\begin{align}
    F(\pi)  = \Phi(x(\mu^\pi,\pi))  - \sum_{t=0}^{T-1}\gamma^{t+1}\sum_{s'\in \setS} \nu^\pi_{t+1,s'} \left[\mu^\pi_{t+1}(s') - \sum_{s\in \setS}\sum_{a\in\setA} P(s'|s, a, \mu^\pi_t, \pi_t)\pi_t(s,a) \mu^\pi_t(s)\right],
\end{align}
Here, we use the original definition of state transition $P$ to distinguish decision variable $\pi$ and the subproblem solution $\mu^\pi$.
By the envelope theorem~\citep{milgrom2002envelope}, the gradient can be explicitly obtained as

\begin{align}
\nabla F(\pi)_{t,s,a} 
= \gamma^t \mu^\pi_t(s)
\Big[R(s,a,\mu^\pi_t,\pi_t)
+ \gamma \sum_{s'\in\setS} \left(P(s'|s,a,\mu^\pi_t,\pi_t) + \mathbb{E}\!\left[\frac{\partial P(s'|\cdot,\mu^\pi_t,\pi_t)}{\partial x_t(s,a)}\right]\right)\nu^\pi_{t+1,s'}\Big], \label{eq:grad_final}
\end{align} 
where $\nabla F(\pi)_{t,s,a}$ is the gradient with respect to $\pi_t(s,a)$.

\indent As per Eq.~\eqref{eq:grad_final}, the factor \(\gamma^t\mu_t^{\pi}(s)\) could be extremely small at states with little mass. 
The same issue is observed in \cite{guo2025markov} when the policy update is weighted by the state visitation frequency. 
To address this problem, we use a state-dependent adaptive step size. Let $\beta>0$ be the base step size and $\bar\beta>0$ be a fixed maximum step size. Then in each iteration $i$, the adaptive step size is set as:
\begin{align}
\eta_{t,s}^{(i)}
=
\begin{cases}
\displaystyle
\min\left\{
\frac{\beta}{\gamma^t\mu_t^{\pi^{(i)}}(s)},
\bar\beta
\right\},
&\mu_t^{\pi^{(i)}}(s)>0,\\
\bar\beta,
&\mu_t^{\pi^{(i)}}(s)=0,
\end{cases}
\end{align}
and the policy is updated by the projected-gradient policy iteration (PG) algorithm:
\begin{align}\label{eq:pg}
    \pi^{(i+1)}_t(s,\cdot) = \text{Proj}_{\Delta(\setA)}\!\left[\pi^{(i)}_t(s,\cdot) + \eta^{(i)}_{t,s} \nabla F(\pi^{(i)})_{t,s,\cdot}\right],
\end{align}
where $\text{Proj}_{\Delta(\setA)}[\cdot]$ is the projection operator.

The complete solution algorithm is summarized in Algorithm~\ref{alg:pg}, while the following theorem states that, with a mild assumption and proper step sizes, the PG algorithm guarantees monotonic improvement of the objective function and every limit point of the policy sequence satisfies the KKT condition of Problem~\eqref{eq:optm_MPPG}. The complete proof is detailed in Appendix~\ref{app:prop_convergence}.

\begin{theorem}\label{prop:convergence}
Suppose the gradient of the objective function $F$ is $L_F$-Lipschitz
continuous on $\Pi$; that is, there exists a constant $L_F>0$ such that
\begin{align}
\left\|
\nabla F(\pi)-\nabla F(\hat{\pi})
\right\|_2
\leq
L_F\left\|
\pi-\hat{\pi}
\right\|_2,
\qquad
\forall \pi,\hat{\pi}\in\Pi.
\end{align}
Let $\beta>0$ and choose
$0<\bar\beta<\frac{2}{L_F}$. 
Then the sequence $\{\pi^{(i)}\}_{i\geq 0}$ generated by the
PG algorithm satisfies
$F(\pi^{(i+1)})> F(\pi^{(i)})$, whenever $\pi^{(i+1)}\neq\pi^{(i)}$.
Moreover, the sequence
$\{\pi^{(i)}\}_{i\geq 0}$ has at least one limit point, and every limit point $\pi^*$ yields  
$(\mu^{\pi^*},\pi^*)$ satisfying the KKT conditions
of Problem~\eqref{eq:optm_MPPG}.
\end{theorem}

\begin{remark}
    It is worth noting that Theorem~\ref{prop:convergence} does not guarantee solution uniqueness. When Problem~\eqref{eq:optm_MPPG} admits multiple KKT points, the PG algorithm may converge to different KKT points depending on the initial policy.
\end{remark}

\begin{algorithm}[h]
\caption{Projected-gradient policy iteration algorithm}
\label{alg:pg}
\begin{algorithmic}

\State \textbf{Inputs:} initial policy $\pi^{(0)}$; initial state distribution $\mu_0$; transition function $P$; reward function $R$; discount factor $\gamma$; maximum iteration $I$; gap tolerance $\epsilon$; base stepsize $\beta$; maximum stepsize $\bar \beta$; initial gap $g =\infty$ and iterate $i= 0$

\State \textbf{Outputs:} equilibrium policy $\pi^*$, mean distribution trajectory $\mu^*$.

\While{$g > \epsilon$ and $i < I$}
\State Solve $\mu^\pi, \nu^\pi$ from KKT conditions at $\pi^{(i)}$. 

\State Evaluate policy gradient as per Eq.~\eqref{eq:grad_final}. 

\State Update policy as per Eq.~\eqref{eq:pg}

\State Calculate gap as $g = \frac{\|\pi^{(i+1)}-\pi^{(i)}\|_1}{\|\pi^{(i)}\|_1}$.

\State Update iterate $i\leftarrow i+1$

\EndWhile

\State Return solutions $\pi^* = \pi^{(i)}$, $\mu^* = \mu^\pi$.

\end{algorithmic}
\end{algorithm}

\section{Connection to other Markovian traffic equilibrium models}\label{sec:connect_MTE}

\subsection{Model connections}
\label{sec:connection_MTE_models}

Similar to the basic idea of MPPG, several existing MTE models admit equivalent optimization formulations~\citep{akamatsu1997decomposition,baillon2008markovian,oyama2022markovian}. A key requirement is the integrability of the link-cost mapping, which is satisfied by 
separable link costs or non-separable costs with a symmetric Jacobian of link costs~\citep{baillon2008markovian,oyama2022markovian}, that is:
\begin{equation}
\frac{\partial c_\ell(f)}{\partial f_m}
=
\frac{\partial c_m(f)}{\partial f_\ell},
\qquad \forall \ell,m,
\end{equation}
where $f$ is the link-flow vector. This parallels Assumption~\ref{ass:sym_reward}, which ensures the existence of a stage potential in MPPG.

An important distinction concerns the dependence of the transition kernel on aggregate flows. 
In conventional link-choice models, transitions are determined by the selected link  and are therefore deterministic and exogenous.
As established in Corollary~\ref{coro:fully-independent}, such transitions satisfy Assumption~\ref{ass:transition}, so a global potential can be constructed as per Theorem~\ref{prop:suffic_potential}.
In contrast, the MPPG framework is more general, as it allows flow-dependent endogenous transition kernels, so long as Assumption~\ref{ass:transition} holds.

Recently, \citet{yao2026perturbed} develop the perturbed utility
Markovian equilibrium (PUME) model, which generalizes the choice model through a convex perturbation function of choice probability, and allows  
both symmetric and asymmetric cost
interactions. PUME is a generalization of MTE, in the following, we therefore establish the connection between PUME and MPPG.
Its exogenous and flow-independent transition kernel satisfies Assumption~\ref{ass:transition}. 
For the finite-horizon formulation, the case of separable link utilities admits the following potential representation, as proved in Appendix~\ref{app:PUME}.
\begin{proposition}\label{prop:PUME}
    Suppose the link utility function $r$ is continuously differentiable and separable, and the perturbation functions $H_s^*$ are continuously differentiable. 
    Then, PUME admits an MPPG representation with stage potential 
    \begin{align}\label{eq:PUME_stage_potential}
        \phi(x_t) = \sum_s\sum_a\int_0^{x_t(s,a)} r(s,a,u) \mathrm{d} u - \sum_s \mu_t(s) H^*_s(\pi_t(s,\cdot)).
    \end{align}
\end{proposition}

\subsection{Algorithm connections}
\label{sec:connection_MTE_algorithms}

Solution algorithms for MTE models mentioned in Section~\ref{sec:connection_MTE_models} typically alternate between network loading and cost updates. For network loading, value iteration is commonly used in the infinite-horizon setting~\citep[e.g.,][]{baillon2008markovian,oyama2022markovian}, while modified policy iteration (MPI) is applied in~\cite{yao2026perturbed} and shows better convergence. As for the cost updates, when an equivalent optimization formulation is available, first-order methods are often applied with accelerations. 
For instance, \citet{baillon2008markovian} proposed gradient-based iterations on either link cost or flow, though reporting slow convergence in both cases. Consequently, they switched to simplified MSA-based link-flow iterations, combined with Newton acceleration.
Differently, \citet{oyama2022markovian} adopted projected gradient for link cost update, combined with momentum acceleration.
As PUME handles more general asymmetric cost interactions, there is no equivalent optimization problem and the equilibrium is formulated as a variational inequality (VI). A meta-algorithm that combines the first-order VI method and safeguarded acceleration is then developed~\citep{yao2026perturbed}.

Since MPPG is defined on a finite horizon, the value functions can be more easily solved from a linear system constructed by a forward network loading. Then, the policy is updated following the projected gradient method.

\section{Numerical experiments}\label{sec:app}

This section presents numerical experiments to validate the
theoretical properties of the proposed MPPG framework and evaluate the performance of the solution algorithm. We construct a
parametric Markov population game that enables systematic control of the endogenous dependence of state transitions and the violation of Assumption~\ref{ass:transition}. We first examine how the transition properties affect
the exact $\varepsilon$ and the theoretical bounds (Section~\ref{subsec:numerical:bounds}), and then compare the proposed algorithm against MSA-based and greedy policy iteration methods, across various state and action space sizes and different time horizons (Section~\ref{subsec:numerical:algorithms}). 

\subsection{A parametric Markov population game}
\label{subsec:numerical:model}

Consider a Markov population game with state space
$\setS=\{1,\ldots,n\}$ and action space
$\setA=\{1,\ldots,m\}$, where $n\geq2$ is even and $m\geq2$. Throughout the numerical experiments, we use the uniform
initial state distribution $\mu_0(s)=1/n$ and set initial policy as $\pi_t^{(0)}(s,a)=1/m, \forall t,s,a$.

Define the immediate and expected reward functions as
\begin{align}
r(s'|s,a,x_t)=R(s,a,x_t)
&= \frac{sa}{nm} -c\left[\frac{x_t(s,a)}{\sum_{a'} x_t(s,a')}
-\frac{\sum_{a'} x_t(s,a')^2}{2(\sum_{a'} x_t(s,a'))^2}\right],\;\forall s'\in\setS,
\label{eq:num:reward}
\end{align}
where $c>0$ is the cost factor associated with state-action flow $x_t(s,a)$. 

It is easy to verify that the expected reward \eqref{eq:num:reward} satisfy Assumption~\ref{ass:sym_reward}, and 
the corresponding stage potential function is 
\begin{align}
    \varphi(x_t)
&=
\sum_{s,a}\frac{sa}{mn}x_t(s,a)
-\frac {c}{2}\sum_s
\frac{\sum_a x_t(s,a)^2}{\sum_a x_t(s,a)}.
\label{eq:num:stagepotential}
\end{align}

We then construct the transition kernel as
\begin{equation}
P(s'|s,a,x_t)=
\frac{1}{n}
\left[
1+z_{s'}q_a(y_t)
\right],
\label{eq:num:transition}
\end{equation}
where $z_s$ is the state-specific indicator defined as $z_s=(-1)^s,\forall s$, which ensures $\sum_{s'} P(s'|s,a,x_t) = 1 + (q_a(y_t)/n)\sum_{s'}z_{s'}=1$. Together with the uniform initial distribution and the range of parameters in $q_a(y_t)$ stated below, it ensures that each state is visited and, accordingly, the reward and stage potential are well defined.

$q_a(\cdot)$ is a function of the vector of aggregate action flows $y_t=(y_{t,a})_{\forall a}$ with
$y_{t,a}=\sum_{s}x_t(s,a)$, $\sum_{a}y_{t,a}=1$ and specified as 
\begin{align}\label{eq:num:action-component}
q_a(y_t)=\kappa_1
\left(\frac{2y_{t,a}}{\|y_t\|}-1\right)
+\kappa_2\left[\frac{\sum_{a_1=1}^{m-1}\sum_{a_2=a_1+1}^{m}
(y_{t,a_1}-y_{t,a_2})^2}{m-1}\right],
\end{align}
where $\kappa_1,\kappa_2\geq 0$ and $\kappa_1+\kappa_2<1$. 

The two components in Eq.~\eqref{eq:num:action-component} are designed to control the satisfaction of Assumption~\ref{ass:transition}. Specifically, the assumption holds when $\kappa_2=0$ and the violation increases with $\kappa_2$. This structural property and the range of $\kappa_1,\kappa_2$ are discussed in Appendix~\ref{app:trans_kernel}.

With the transition kernel defined above, a uniform upper bound of the Dobrushin contraction coefficient can be derived as
\begin{align}
\theta_t(x_t)&=\frac{1}{2}\max_{\substack{s_1,s_2\in\setS\\a_1,a_2\in\setA}}\sum_{s'\in\setS}\left|P(s'\mid s_1,a_1,x_t)-P(s'\mid s_2,a_2,x_t)\right|=\frac{1}{2}\max_{a_1, a_2\in\mathcal A} \sum_{s'\in\setS}\left|\frac{z_{s'}}{n}q_{a_1}(y_t)-\frac{z_{s'}}{n}q_{a_2}(y_t)\right| \nonumber\\
&=\kappa_1 \max_{a_1, a_2\in\mathcal A} \left|\frac{y_{t,a_1}}{\|y_t\|} - \frac{y_{t,a_2}}{\|y_t\|}\right| = 
\kappa_1\frac{\max_a y_{t,a}-\min_a y_{t,a}}{\|y_t\|_2}\leq \kappa_1,
\label{eq:num:dobrushin}
\end{align}
which will be used in the numerical experiments when evaluating the theoretical bounds.

Similarly, we construct the reward bounds as follows:
\begin{align}
   R(s,a,x_t) &= \frac{sa}{nm} - c\left[\frac{x_t(s,a)}{\sum_{a'} x_t(s,a')}
-\frac{\sum_{a'} x_t(s,a')^2}{2(\sum_{a'} x_t(s,a'))^2}\right] \nonumber\\
&= \frac{sa}{nm} - c\left[\pi_t(s,a) - \frac{1}{2}\sum_{a'}\pi_t(s,a')^2\right] \leq  \frac{sa}{nm} +\frac{c}{2}\sum_{a'}\pi_t(s,a')^2 \leq 1+\frac{c}{2}\\
& = \frac{sa}{nm} +\frac{c}{2}\left[(\pi_t(s,a)-1)^2 + \sum_{a'\neq a} \pi_t(s,a')^2 - 1 \right] \geq \frac{sa}{nm} - \frac{c}{2} \geq \frac{1}{nm} - \frac{c}{2},
\end{align}
and the reward span is $\Delta r = 1 - \frac{1}{nm} + c$.

\subsection{Validation of theoretical results}
\label{subsec:numerical:bounds}
We evaluate the $\varepsilon$-equilibrium and approximation bounds established
in Propositions~\ref{theo:priori bound} and~\ref{theo:posteriori_bound} by varying the transition parameters
$\kappa_1$ and $\kappa_2$. We consider an instance with
$n=6$, $m=6$, and $T=10$.
The baseline parameter values and their respective
test ranges are summarized in Table~\ref{tab:num_parameters}.
For each parameter configuration, we apply the proposed
PG algorithm and evaluate the realized equilibrium gap
$\varepsilon$, together with its a priori bound $\varepsilon_{\mathrm{pri}}$ and a posteriori
bound $\varepsilon_{\mathrm{post}}$. 
Besides, since we focus on the effect of transition dynamics on the approximated equilibrium, we keep the discount factor $\gamma$ and reward parameter $c$ fixed throughout the experiments. 

\begin{table}[htbp]
\centering
\caption{Parameters baseline and test ranges.}
\label{tab:num_parameters}
\begin{tabular}{lcc}
\hline
Parameter & Baseline & Tested values \\
\hline
$\kappa_1$ & $0.2$ & $0,\,0.1,\,0.2,\,0.3,\,0.4,\,0.5,\,0.6,\,0.7$ \\
$\kappa_2$ & $0.2$ & $0,\,0.1,\,0.2,\,0.3,\,0.4,\,0.5,\,0.6,\,0.7$ \\
$\gamma$ & $1.0$ & -- \\
$c$ & $1.0$ & --  \\
\hline
\end{tabular}
\end{table}

\begin{figure}[h]
    \centering
    \begin{subfigure}[t]{0.46\textwidth}
        \centering
        \includegraphics[width=\textwidth]{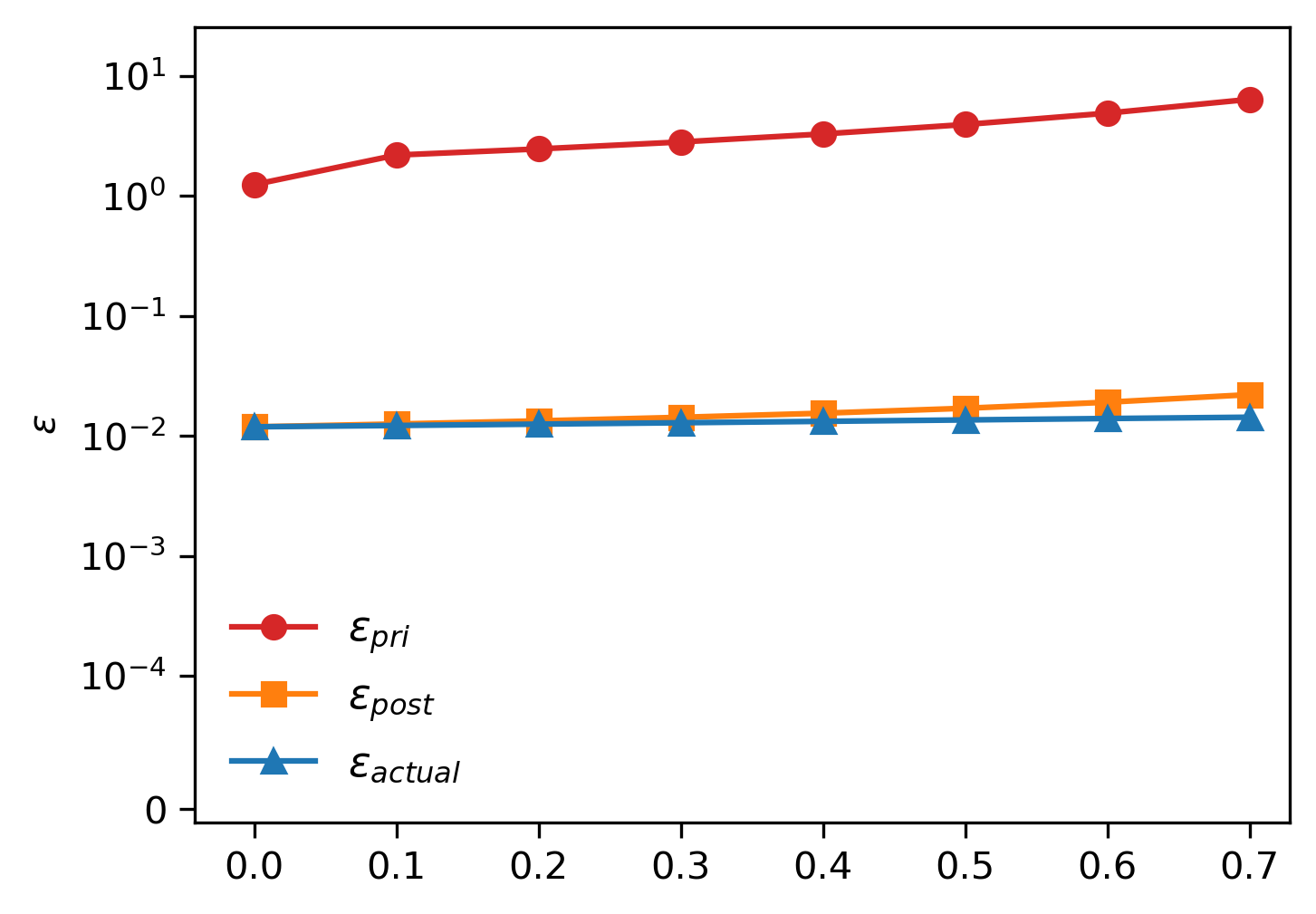}
    \caption{$\kappa_1$}\label{fig:kappa1_sweep}
    \end{subfigure}
    \begin{subfigure}[t]{0.46\textwidth}
        \centering
    \includegraphics[width=\textwidth]{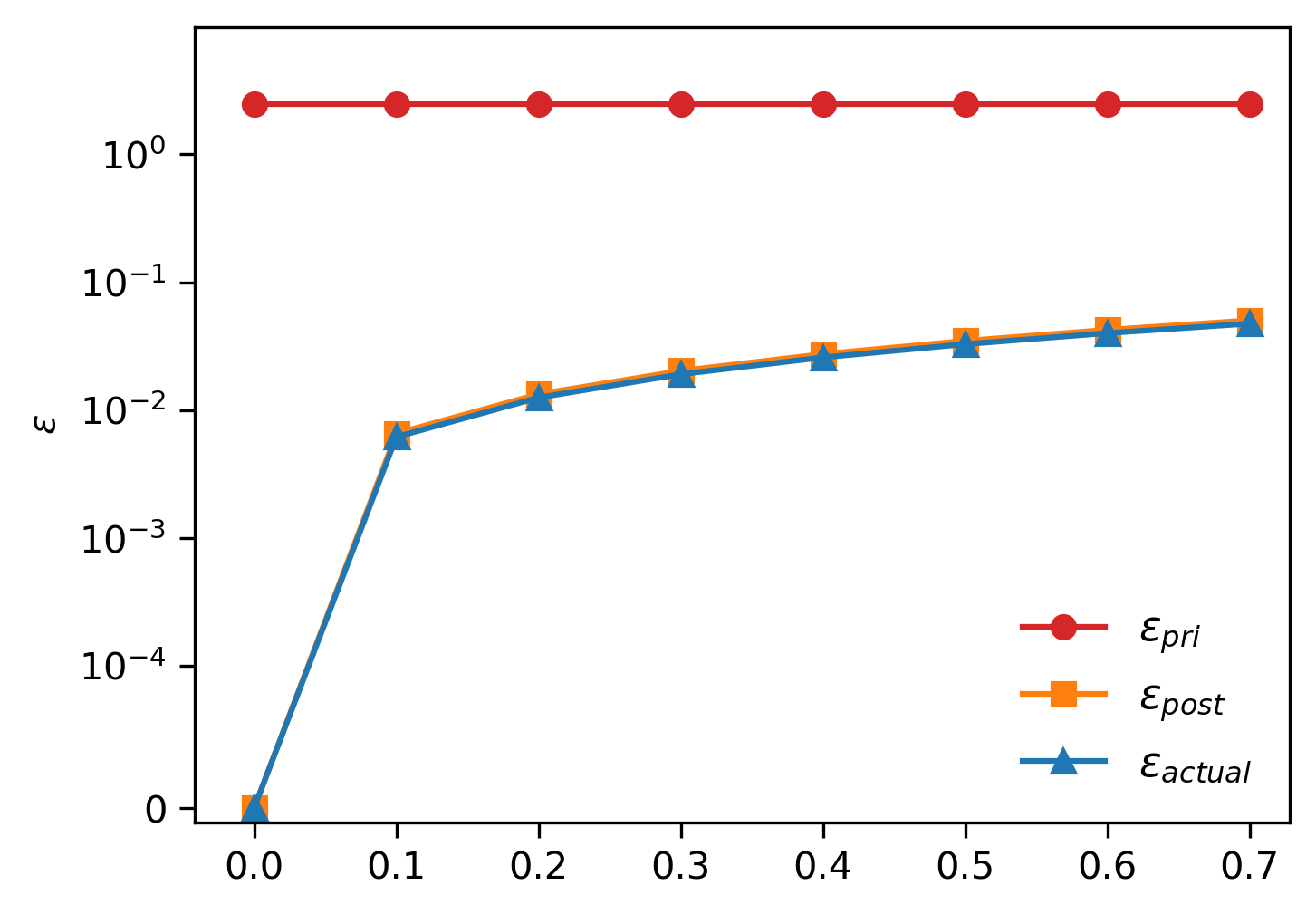}
        \caption{$\kappa_2$}\label{fig:kappa2_sweep}
    \end{subfigure}
    \caption{$\varepsilon$ and theoretical bounds under various parameters }\label{fig:bound_para}
\end{figure}

The results are presented in Fig.~\ref{fig:bound_para}. Across all parameter configurations, the a priori bound is substantially more conservative than the a posteriori bound. As shown in Fig.~\ref{fig:kappa1_sweep}, increasing $\kappa_1$ leaves the realized $\varepsilon$ nearly unchanged, as it does not affect the violation of Assumption~\ref{ass:transition}. However, it increases the Dobrushin coefficient $\theta$, resulting in larger theoretical bounds.

In contrast, $\kappa_2$ directly affects the violation of Assumption~\ref{ass:transition}. 
Specifically, the assumption holds when $\kappa_2=0$, thus both $\varepsilon$ and the a posteriori bound are zero, as shown in Fig.~\ref{fig:kappa2_sweep}. 
Nevertheless, the a priori bound remains relatively large due to 
the first term in Eq.~\eqref{eq:epsilon_apriori_cases}. 
When $\kappa_2>0$, $\varepsilon$ jumps from zero to a positive value, though it does not further increase substantially as $\kappa_2$ keeps growing, demonstrating robustness of the approximation error.

\subsection{Algorithmic convergence and scalability}
\label{subsec:numerical:algorithms}

We further evaluate the computational performance of the
proposed projected-gradient policy iteration (PG) by comparing
it with MSA-based policy iteration (MSA) and greedy policy
iteration (Greedy). The three algorithms are tested on the
same parametric Markov population game under varying time
horizons, state-space sizes, and action-space sizes.

The PG algorithm follows the procedure developed in
Section~\ref{sec:mppg_algo}. We fix the base step size as $\beta=1$, and test several values of the maximum step size $\bar \beta$. Since the qualitative conclusions of the algorithm comparison remain consistent across these step size settings, we only report results of one setting. Both MSA and Greedy include the same forward propagation as PG, but employ different policy update rules.

Let $\widehat{\pi}^{(k)}$ denote the best response to the
mean field induced by the current policy $\pi^{(k)}$. MSA and Greedy methods update policy as follows:
\begin{align}
\text{MSA: }\quad & \pi^{(k+1)}= (1-\alpha_k)\pi^{(k)} +\alpha_k\widehat{\pi}^{(k)}, \\
\text{Greedy: }\quad & \pi^{(k+1)}=\widehat{\pi}^{(k)},
\label{eq:num:policy_iteration}
\end{align}
where $\alpha_k=\frac{1}{1+k}$ is the diminishing step size.
Based on the results in Fig.~\ref{fig:bound_para}, all parameters are set to their baseline values in Table~\ref{tab:num_parameters}, under which the realized $\varepsilon$ and its bounds are relatively stable.

To assess convergence to equilibrium, besides the policy gap, we also consider the normalized Q-gap
\begin{equation}
g_{Q}(\pi)=
\frac{
\displaystyle
\sum_{t\in\setT}
\sum_{s\in\setS}
\sum_{a\in\setA_+(t,s;\pi)}
\left|
Q_t(s,a,\mu,\pi)
-
\max_{a'\in\setA_+(t,s;\pi)}
Q_t(s,a',\mu,\pi)
\right|
}{
\displaystyle
\sum_{t\in\setT}
\sum_{s\in\setS}
\sum_{a\in\setA_+(t,s;\pi)}
\left|
Q_t(s,a,\mu,\pi)
\right|},
\label{eq:Q_gap}
\end{equation}
which measures the total difference between Q-values of active actions and the maximum Q-value at the same state, normalized by the total absolute Q-value. 
We note that true Q-values are used to compute the Q-gap for PG, although the policy update is based on pseudo Q-values.

\begin{figure}[htb]
    \centering
    \includegraphics[width=\linewidth]{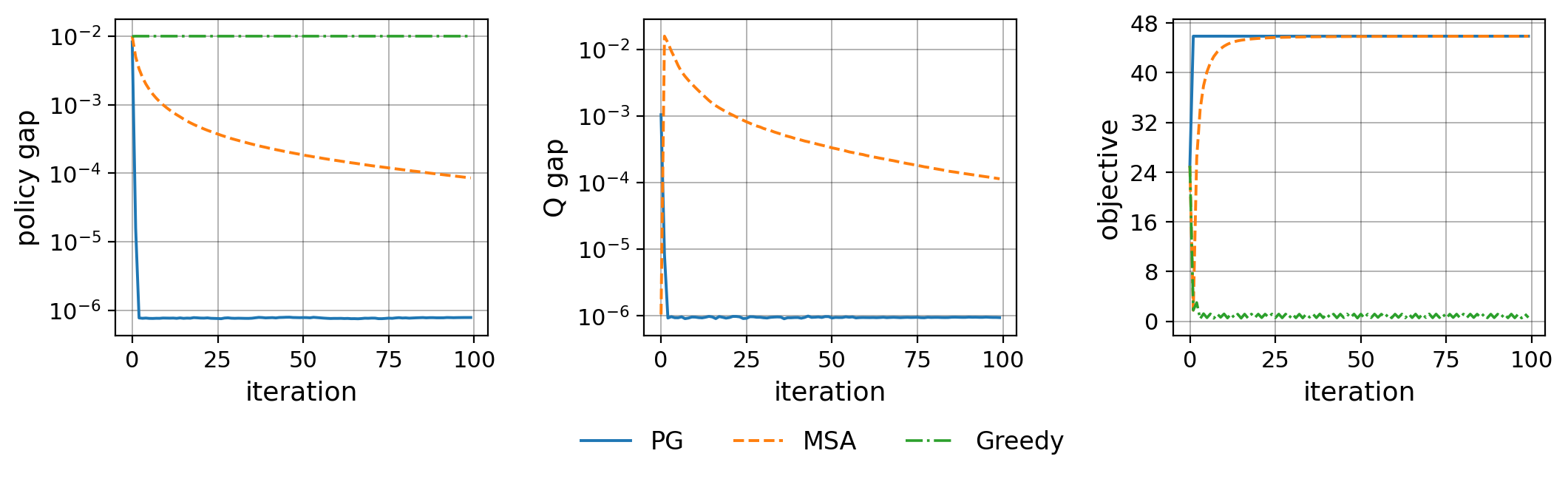}
    \caption{Convergence of three algorithms, with $T=100$, $n\times m=200\times 200$}
    \label{fig:conv_scale_example}
\end{figure}

Figure~\ref{fig:conv_scale_example} compares the convergence behavior of PG,
MSA, and Greedy methods, tracking the
policy gap, Q-gap, and objective value defined
in~\eqref{eq:potential_construct} over iterations. Due to space constraints, we focus on the case $T=100$ and
$n\times m=200\times200$, while the full results are
provided in Appendix~\ref{app:numerical}. 
The Q-gap of Greedy is omitted because it is
always zero due to the greedy update. Both PG and MSA exhibit convergence, whereas Greedy fails to converge with an almost constant policy gap over iterations. Besides, the objective value under Greedy oscillates at a lower level, whereas PG and MSA quickly converge to a nearly identical objective value after a few iterations, providing empirical
support for the proposed objective formulation.

Figure~\ref{fig:conv_scale_example} also demonstrates a clear advantage of PG over MSA in reducing both the policy
gap and Q-gap. During the initial phase, both gaps under PG decrease approximately exponentially (linearly on a logarithmic scale), falling below $10^{-6}$, within 10 iterations.
In contrast, the decay of gaps under MSA is substantially slower, and both gaps remain considerably larger than those of PG at the end of the iterations.

To assess whether this convergence advantage persists across
problem settings, we test on different problem scales and report the iterations required to achieve a Q-gap below $10^{-3}$ in Table~\ref{tab:pg-msa-comparison}. For MSA, we check the $Q$-gap from the second update,
because the first update  produces a deterministic policy
with zero $Q$-gap. It can be found that PG consistently requires substantially fewer iterations than MSA, while MSA fails to reach this tolerance within 100 iterations in some settings with small $T$. It is worth noting that the iterations decrease with $T$ because, as cumulative rewards, Q-values increase with time horizon, leading to a smaller normalized gap at the same absolute difference.

\begin{table}[htbp]
    \centering
    \begin{threeparttable}
    \caption{Iterations required to reach a $Q$-gap below $10^{-3}$.}
    \label{tab:pg-msa-comparison}
    \setlength{\tabcolsep}{7pt}
    \renewcommand{\arraystretch}{1.15}
    \begin{tabular}{@{}c *{4}{rr}@{}}
        \toprule
        & \multicolumn{8}{c}{$n\times m$} \\
        \cmidrule(lr){2-9}
        $T$
        & \multicolumn{2}{c}{$30\times30$}
        & \multicolumn{2}{c}{$50\times50$}
        & \multicolumn{2}{c}{$100\times100$}
        & \multicolumn{2}{c}{$200\times200$} \\
        \cmidrule(lr){2-3}
        \cmidrule(lr){4-5}
        \cmidrule(lr){6-7}
        \cmidrule(lr){8-9}
        & PG & MSA & PG & MSA & PG & MSA & PG & MSA \\
        \midrule
         20 & 2 & $>100$ & 2 & $>100$
            & 2 & $>100$ & 2 & $>100$ \\
         50 & 2 & $47$ & 2 & 46
            & 2 & 46       & 2 & 43 \\
         80 & 2 & 31       & 2 & 31
            & 2 & 29       & 2 & 28 \\
        100 & 2 & 26       & 2 & 25
            & 2 & 25       & 2 & 23 \\
        \bottomrule
    \end{tabular}
    
    \end{threeparttable}
\end{table}

\begin{figure}[h]
    \centering
    \begin{subfigure}[t]{0.49\textwidth}
        \centering
        \includegraphics[width=\textwidth]{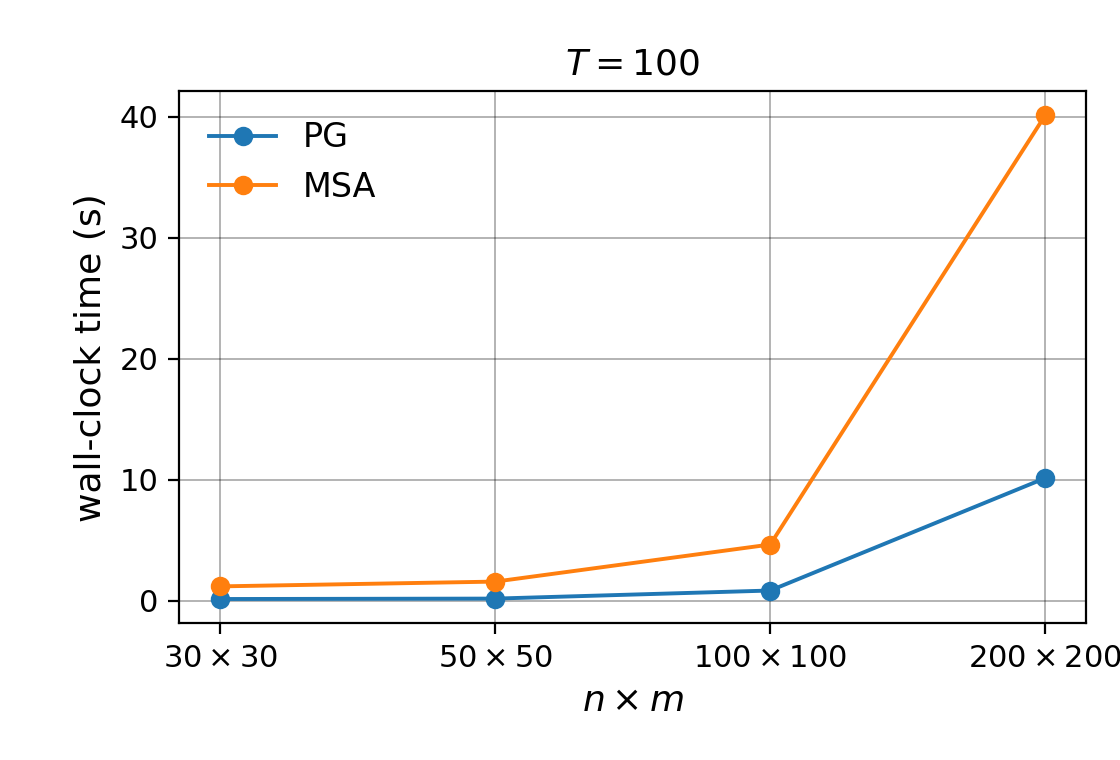}
    \caption{varying $n\times m$}\label{fig:wc_time_SA}
    \end{subfigure}
    \begin{subfigure}[t]{0.49\textwidth}
        \centering
    \includegraphics[width=\textwidth]{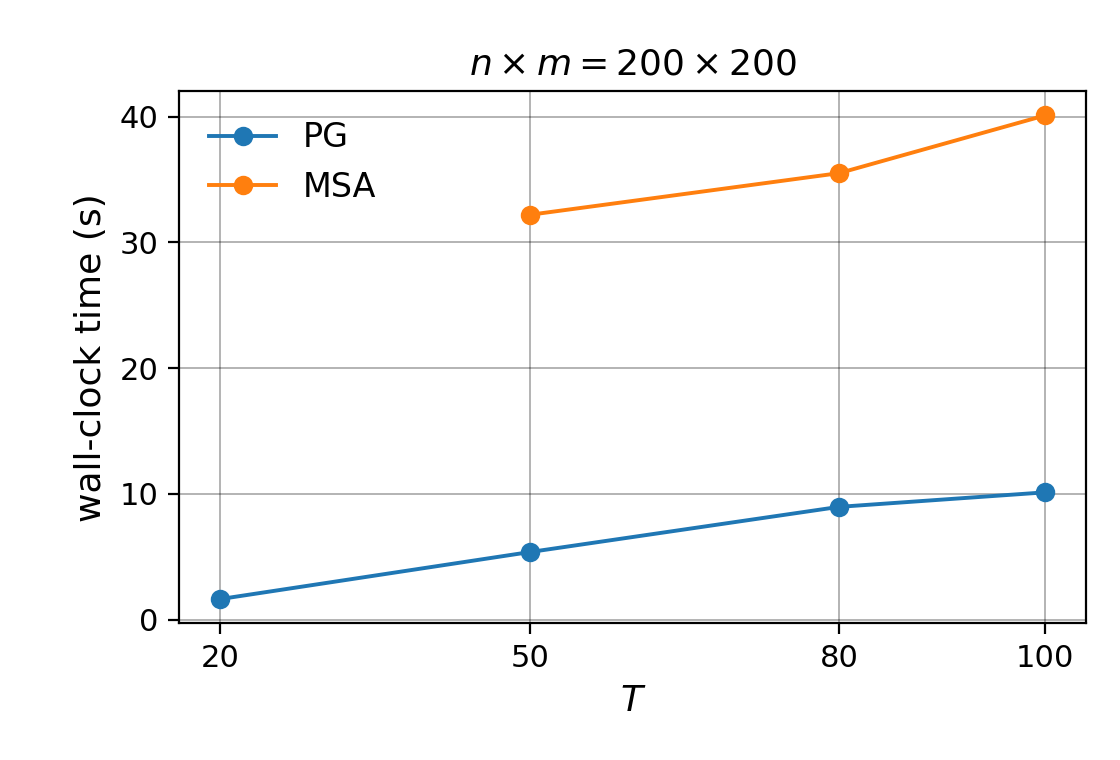}
        \caption{varying $T$}\label{fig:wc_time_T}
    \end{subfigure}
    \caption{Wall-clock time under different $n\times m$ and $T$}\label{fig:wall clock time}
\end{figure}

Figure~\ref{fig:wall clock time} further compares the wall-clock time required by PG and MSA to achieve the Q-gap threshold in Table~\ref{tab:pg-msa-comparison} as either $T$ or the state-action space increases.
As expected, the computation time increases as the problem scale expands, and the jump is particularly significant when the state-action space grows up to $200\times200$. On the other hand, a longer time horizon only leads to a mild increase in computation time. 
Overall, PG consistently outperforms MSA in both sensitivity analyses. 
Specifically, Fig.~\ref{fig:wc_time_T} shows that, at $n\times m=200\times200$ and $T=20$, MSA fails to reach the desired gap threshold within 100 iterations. 
These results demonstrate the computational efficiency of PG and provide empirical evidence of its scalability over the tested range.

\section{Conclusions}\label{sec:conclusion}
This study introduces the Markov population potential game (MPPG) as a unified framework for modeling large-scale, dynamic, and strategic interactions among many self-interested agents. We establish a set of sufficient conditions under which a Markov population game is an MPPG, which better balances generality and ease of verification compared to existing ones in the literature. 
These conditions also allow us to explicitly construct the potential function and prove that each KKT point of the associated optimization problem corresponds to the equilibrium of MPPG. 
We further show that if the condition on state transitions does not hold, such a KKT point yields a support-restricted $\varepsilon$-equilibrium and provide the a priori and posteriori bounds of $\varepsilon$. 
On the computational side, we develop a  projected-gradient policy iteration (PG) algorithm.

The MPPG framework is validated through numerical experiments with systematic control over the degree of condition violation. 
By varying the model parameters, we examine how the transition kernel affects the realized equilibrium approximation error $\varepsilon$ and assess the tightness of the a priori and a posteriori bounds. We also compare PG, MSA, and greedy-based policy iteration methods across different state–action space sizes and time horizons. The results demonstrate the computational efficiency and scalability of PG, which is consistently faster and more stable in convergence across the tested configurations.

Besides the connections to classical Markovian traffic equilibrium (MTE) models and the recent perturbed utility Markovian equilibrium (PUME) are discussed in this paper. 
The MPPG framework is also closely related to several general and application-specific models proposed in the literature. For instance, the user equilibrium with recourse~\citep{unnikrishnan2009user} also follows the MTE framework with flow-dependent link costs and deterministic state transitions. 
Differently,  the parking search equilibrium model developed in~\cite{boyles2015parking} considers fixed link costs but endogenous flow-dependent parking-success probabilities. It thus naturally satisfies Assumption~\ref{ass:sym_reward} but not necessarily Assumption~\ref{ass:transition}. 
Similarly, RIVER~\citep{zhang2023ride} and eRIVER~\citep{zhang2023routing}, developed to model ride-hailing operations, also satisfy Assumption~\ref{ass:sym_reward} with a symmetric Jacobian of the reward function but not Assumption~\ref{ass:transition}. In these cases, an optimization can still be constructed based on the stage potential, though its solution only corresponds to an $\varepsilon$-equilibrium. 
In contrast, CARMA~\citep{elokda2025carma}, a recently proposed non-tradable credit mechanism and also built on a Markov population game framework, satisfies neither Assumption~\ref{ass:sym_reward} nor Assumption~\ref{ass:transition}. Therefore, it cannot directly apply the analysis and solution algorithm presented in this study.

There are several directions for future research. First, taking advantage of the potential, second-order solution algorithms can be developed to further enhance the scalability of MPPG. Second, the closed-form potential characterization derived in our model offers a solid theoretical foundation for exploring incentive and mechanism design in large-scale transportation systems. Moreover, MPPG is closely related to recent studies on $\alpha$-MPG \citep{guo2025markov} and more general dynamic games \citep{Guo2025Alpha}. It is thus interesting to explore the relationship between MPPG and these emerging game-theoretical frameworks.

\Urlmuskip=0mu plus 1mu\relax
\bibliographystyle{abbrvnat}
\bibliography{ref_mpg,ref}

\newpage
\appendix
\section*{Appendix}

\section{Proof of Theorem~ \ref{prop:existence_equ}}\label{app:existence_equ}
Define the joint space of policy and state distribution as $X := \mathcal M \times \Pi$, then the equilibrium conditions can be written as a fixed-point problem on $X$ by defining a correspondence $F:X\rightrightarrows X$:
\begin{align}
    F(\mu,\pi) :=\Bigl\{(\hat{\mu},\hat{\pi})\in X \;\Big|\;
    \begin{array}{l}
    \hat{\pi}_t(s, \cdot)
    \in \arg\max_{\pi'_t(s, \cdot)\in \Delta(\mathcal{A})}\sum_{a\in\setA}\pi'_t(s,a)Q_t(s,a,\mu,\pi),
    \quad \forall t,s,\\
    \hat{\mu}_{t+1}(s')
    = \sum_{s\in \setS}\sum_{a\in \setA}
    P(s'|s,a,\mu_t,\pi_t)\mu_t(s)\pi_t(s,a),
    \quad \forall t,s'
    \end{array}
    \Bigr\},
\end{align}
Accordingly, $(\mu^*,\pi^*)$ is an equilibrium if and only if
$(\mu^*,\pi^*)\in F(\mu^*,\pi^*)$. Then, proving the existence of equilibrium is equivalent to proving the existence of a fixed point. To this end, we invoke Kakutani’s fixed-point theorem by showing that i) $X$ is nonempty, compact, and convex, and ii) $F$ is upper hemicontinuous, nonempty, closed and convex on $X$. 

The first condition is easy to prove since $\Delta(\setA)$ and $\Delta(\setS)$ are simplices over finite sets, and these properties are preserved under finite Cartesian products.
To prove the second condition, we first note that the state dynamics are continuous and given the initial distribution, $\hat{\mu}$ is uniquely determined by any $(\mu,\pi) \in X$. It is easy to see the correspondence of $\mu$ is upper hemicontinuous (indeed, continuous), non-empty, compact, and convex valued on $X$.
We then prove the set of optimal policy $\hat{\pi}$ is also nonempty, closed and convex and upper hemicontinuous, which then induces the same property of $F$. 
Due to the continuity of $R$ and $P$, the function $Q_t(\cdot,\cdot;\mu,\pi)$ is continuous in $(\mu,\pi)$, which can be shown by a standard backward induction over finite horizons. 

Given that the set of feasible policies $\Delta(\setA)$ is nonempty, compact and convex, the optimal policy $\hat{\pi}$ always exists. And we have by Berge's maximum theorem that  the best-response correspondence of $\hat{\pi}$ is upper hemicontinuous and compact (thus closed).
Besides, it is easy to show that the convex combination of any two optimal policies is also optimal. Hence, the set of $\hat{\pi}$ is also convex.

Finally, using the definition, we have that $F(\mu, \pi)$, a finite Cartesian product of nonempty, compact, convex, and upper hemicontinuous correspondences, is nonempty, compact-valued, convex-valued, and upper hemicontinuous.

\section{Proof of Theorem~\ref{prop:suffic_potential}}\label{app:prop_suffic}

Assumption \ref{ass:sym_reward} states that $R$ is continuously differentiable and symmetric in the marginal effect. This is exactly the \textit{externality symmetry} condition in population potential games ~\citep{sandholm2001potential, sandholm2009large}, which
leads to the existence of a stage potential function $\phi(x_t)$ satisfying Eq.~\eqref{eq:stage_potential}. Then we need to prove that, with Assumption \ref{ass:transition}, a global potential $\Phi$ constructed as
\begin{equation}
    \Phi(x) = \sum_{t=0}^{T-1} \gamma^t \phi(x_t),
\end{equation}
satisfies Condition~\eqref{eq:MPPG}.

We prove this result by backward induction. 
At the terminal stage $t = T -1$, the derivative of the global potential $\Phi$ with respect to $x_{T-1}(s,a)$ is given by
\begin{align}
    \frac{\partial \Phi(x)}{\partial x_{T-1}(s,a)} 
    = \gamma^{T-1} \frac{\partial \phi(x_{T-1})}{\partial x_{T-1}(s,a)} + \sum_{\tau=0}^{T-2}\gamma^{\tau} \frac{\partial \phi(x_{\tau})}{\partial x_{T-1}(s,a)} = \gamma^{T-1} R(s,a,x_{T-1}),
\end{align}
where the equality holds because early stage potential is independent of future state and policy and the stage potential satisfies Eq.~\eqref{eq:stage_potential}. 

By convention, the terminal value is normalized as \(V_T(s,x) = 0\) for any \(s \in \mathcal{S}\).
Hence, the Q-value at $t = T-1$ reduces to $Q_{T-1}(s,a,x) = R(s,a,x_{T-1})$.
It thus concludes that $\Phi$ satisfies Condition~\eqref{eq:MPPG} at the terminal stage.

Assume for period $t+1$, it holds that
\begin{equation}
    \frac{\partial \Phi(x)}{\partial x_{t+1}(s',a')}
    = \gamma^{t+1} Q_{t+1}(s',a',x),
    \quad \forall s', a'.
\end{equation}
We now prove the result also holds for period $t$. To start with, we expand the derivative of $\Phi$ with respect to $x_t$ as 
\begin{align}\label{eq:proof-prop1-potential}
    \frac{\partial \Phi(x)}{\partial x_t(s,a)}
    &= \gamma^t\frac{\partial \phi(x_t)}{\partial x_t(s,a)}
    + \sum_{\tau=t+1}^{T-1}\gamma^\tau\frac{\partial \phi(x_\tau)}{\partial x_t(s,a)} + \sum_{\tau'=0}^{t-1}\gamma^{\tau'}\frac{\partial \phi(x_{\tau'})}{\partial x_t(s,a)}, \notag \\
    & = \gamma^t\frac{\partial \phi(x_t)}{\partial x_t(s,a)}
    + \sum_{\tau=t+1}^{T-1}\gamma^\tau\frac{\partial \phi(x_\tau)}{\partial x_t(s,a)},
\end{align}
where the second equality holds because the past stage-potentials are independent of future state and policy.
As per Eq.~\eqref{eq:stage_potential}, the first term reduces to $\gamma^t R(s,a,x_t)$.
The second term is further expanded as 
\begin{equation}\label{eq:proof-prop1-future}
    \sum_{\tau=t+1}^{T-1}\gamma^\tau\frac{\partial \phi(x_\tau)}{\partial x_t(s,a)}
    = \sum_{s'\in\mathcal{S}} \left[\sum_{\tau=t+1}^{T-1}\gamma^\tau\frac{\partial \phi(x_\tau)}{\partial\mu_{t+1}(s')} \right]
      \frac{\partial \mu_{t+1}(s')}{\partial x_t(s,a)},
\end{equation}
and the factor in the bracket is further expressed as
\begin{align}
    \sum_{\tau=t+1}^{T-1}\gamma^\tau\frac{\partial \phi(x_\tau)}{\partial\mu_{t+1}(s')} &= \sum_{a' \in \mathcal{A}} \left[ \sum_{\tau=t+1}^{T-1}\gamma^\tau\frac{\partial \phi(x_\tau)}{\partial x_{t+1}(s',a')} + \sum_{\tau=0}^{t}\gamma^\tau\frac{\partial \phi(x_\tau)}{\partial x_{t+1}(s',a')}\right] \frac{\partial x_{t+1}(s',a')}{\partial \mu_{t+1}(s')} \notag \\
    &= \sum_{a' \in \mathcal{A}} \frac{\partial \Phi(x)}{\partial x_{t+1}(s',a')} \pi_{t+1}(s',a') =  \gamma^{t+1} \sum_{a' \in \mathcal{A}}  Q_{t+1}(s',a',x) \pi_{t+1}(s',a') \notag \\
    &=  \gamma^{t+1} V_{t+1}(s',x),
\end{align}
where the first equality is, again, due to the independence of past-stage potential on future stage and policy, the second equality follows by the definition of $\Phi$, the third equality applies the induction assumption at $t+1$, and the last equality uses the definition of $V$.

Next, we expand the second factor in Eq.~\eqref{eq:proof-prop1-future} as 
\begin{align}\label{eq:proof-prop1-transition}
     \frac{\partial \mu_{t+1}(s')}{\partial x_t(s,a)} &= \frac{\partial}{\partial x_t(s,a)} \left[ \sum_{s''\in\setS}\sum_{a''\in\setA} P(s'|s'', a'', x_t) x_t(s'',a'')\right]\nonumber\\
     &= P(s'|s,a,x_t) + \sum_{s''\in\setS}\sum_{a''\in\setA} \frac{\partial P(s'|s'', a'', x_t)}{\partial x_t(s,a)}x_t(s'',a''). 
\end{align}
Under Assumption~\ref{ass:transition}, the second term in Eq.~\eqref{eq:proof-prop1-transition} is simply zero. Accordingly, Eq.~\eqref{eq:proof-prop1-potential} reduces to 
\begin{align}
    \frac{\partial \Phi(x)}{\partial x_t(s,a)}
    = \gamma^t R(s,a,x_t)
       + \gamma^{t+1} \sum_{s'\in\mathcal{S}}P(s'|s,a,x_t) V_{t+1}(s',x) = \gamma^t Q_t(s,a,x).
\end{align}
The result thus holds for $t$, and the proof is completed by induction.

\section{Proof of Corollaries~\ref{coro:fully-independent} and ~\ref{coro:symm-independent}}\label{app:coro_transition}
It is easy to show that Eq.~\eqref{eq:indep-transition} implies Assumption~\ref{ass:transition}, thus the result holds. 
The condition~\eqref{eq:symm-independent} is expanded as 
\begin{align}
 \nabla P(s'|s,a,x_t)^\top x_t &= \sum_{s''\in\setS}\sum_{a''\in\setA} \frac{\partial P(s'|s,a,x_t)}{\partial x_t(s'',a'')} x_t(s'',a'') = \sum_{s''\in\setS}\sum_{a''\in\setA} \frac{\partial P(s'|s'',a'',x_t)}{\partial x_t(s,a)} x_t(s'',a'')\nonumber\\
&=0
\end{align}
where the second equality is due to symmetric state transition. It thus leads to Assumption~\ref{ass:transition} and completes the proof.

\section{Proof of Theorem~\ref{prop:optm_MPPG}}\label{app:prop_optm_MPPG}

To prove a KKT point of Problem~\eqref{eq:optm_MPPG} also satisfies the equilibrium conditions \eqref{eq:def_equ}, we first  rewrite Problem~\eqref{eq:optm_MPPG} as follows:
\begin{subequations}\label{eq:optm_MPPG_reform}
    \begin{align}
        \max_{\mu,\pi} \quad & \Phi(x(\mu,\pi)),\\
        s.t.\quad  &  \gamma^{t+1}\left(\mu_{t+1}(s') -\sum_{s\in \setS}\sum_{a\in\setA} P(s'|s, a, x_t(\mu, \pi)) \mu_t(s) \pi_t(s,a)\right) = 0,\quad  \forall s'\in\setS, t\in\setT, \label{eq:optm_MPPG_dynm_equiv}\\ 
        & \gamma^t \mu_t(s)\left( \sum_{a\in \setA} \pi_t(s,a)-1\right) = 0 \quad \forall s\in\setS, t\in\setT,\label{eq:optm_MPPG_simplex_equiv}\\
        & \gamma^t \mu_t(s) \pi_t(s,a)\geq 0, \quad \forall s\in\setS, a\in\setA, t\in\setT.\label{eq:optm_MPPG_nonneg_equiv}
    \end{align}
\end{subequations}
Note that Problem~\eqref{eq:optm_MPPG_reform} is equivalent to \eqref{eq:optm_MPPG} only if $\mu_t(s) > 0, \forall t \in \mathcal{T}, s \in \mathcal{S}$. When $\mu_t(s) = 0$, the policy $\pi_t(s,a)$ is not effectively constrained. Nevertheless, it does not cause any issue to the corresponding MPPG because a state with zero population mass has no impact on the game dynamics.  
    
The Lagrangian of Problem~\eqref{eq:optm_MPPG_reform} is then written as 
\begin{align}
    \setL(\mu,\pi,\nu, \lambda, \eta) = & \Phi(x(\mu,\pi)) - \sum_{t=0}^{T-1}\gamma^{t+1}\sum_{s'\in \setS} \nu_{t+1,s'} \left(\mu_{t+1}(s') - \sum_{s\in \setS}\sum_{a\in\setA} P(s'|s, a, x_t(\mu,\pi)) \mu_t(s)\pi_t(s,a)\right)\nonumber\\
    & - \sum_{t=0}^{T-1}\gamma^t\sum_{s\in \setS} \lambda_{t,s} \mu_t(s) \left(\sum_{a\in \setA} \pi_t(s,a) - 1 \right) + \sum_{t=0}^{T-1}\gamma^t\sum_{s\in \setS}\sum_{a\in \setA} \eta_{t,sa} \mu_t(s)\pi_t(s,a),
\end{align}
where $\nu,\lambda,\eta$ are Lagrange multipliers for Constraints~\eqref{eq:optm_MPPG_dynm_equiv}, \eqref{eq:optm_MPPG_simplex_equiv}, and \eqref{eq:optm_MPPG_nonneg_equiv}, respectively.

At KKT points,
the primal and dual variables satisfy the following conditions:
\begin{subequations}\label{eq:KKT}
    \begin{align}
        \frac{\partial \setL}{\partial \pi_t(s,a)} &= \frac{\partial \Phi(x(\mu,\pi))}{\partial \pi_t(s,a)} + \gamma^{t+1}\sum_{s'\in\setS}\nu_{t+1,s'}\sum_{s''\in\setS}\sum_{a''\in\setA} \frac{\partial P(s'|s'', a'', x_t(\mu,\pi)) \mu_t(s'')\pi_t(s'',a'')}{\partial \pi_t(s,a)} \nonumber\\
         &\quad - \gamma^t\lambda_{t,s} \mu_t(s) + \gamma^t \eta_{t,sa}\mu_t(s) = 0,\label{eq:KKT_grad_pi}\\
        \frac{\partial \setL}{\partial \mu_t(s)} &= \frac{\partial \Phi(x(\mu,\pi))}{\partial \mu_t(s)} -  \gamma^t \nu_{t,s} + \gamma^{t+1} \sum_{s'\in\setS} \nu_{t+1,s'} \sum_{s''\in \setS}\sum_{a''\in\setA} \frac{\partial P(s'|s'', a'', x_t(\mu,\pi))\mu_t(s'')\pi_t(s'',a'') }{\partial \mu_t(s)} \nonumber\\
        & \quad -\gamma^t \lambda_{t,s}\left(\sum_{a\in \setA} \pi_t(s,a) - 1 \right) + \gamma^t\sum_{a\in\setA} \eta_{t,sa}\pi_t(s,a) =0 \label{eq:KKT_grad_mu}\\
        \mu_{t+1}(s') &= \sum_{s\in \setS}\sum_{a\in\setA} P(s'|s, a, x_t(\mu,\pi)) \mu_t(s)\pi_t(s,a) ,\label{eq:KKT_eq_mu}\\
        \sum_{a\in \setA} \pi_t(s,a) &= 1,\label{eq:KKT_eq_pi}\\
        \eta_{t,sa}\pi_t(s,a) &= 0,\label{eq:KKT_comple}\\
        \eta_{t,sa} & \geq 0.\label{eq:KKT_nonneg}
    \end{align}
\end{subequations}
Eq.~\eqref{eq:KKT_eq_mu} is exactly the equilibrium condition about state distribution dynamics.
Eqs.~\eqref{eq:KKT_grad_pi} and \eqref{eq:KKT_grad_mu} are further expanded and rearranged as 
\begin{subequations}
    \begin{align}
        0 &=   \gamma^t\mu_t(s) R(s,a,x_t(\mu,\pi)) - \gamma^t\lambda_{t,s} \mu_t(s) +  \gamma^t\eta_{t,sa} \mu_t(s) \label{eq:simplified_KKT_lambda} \\
        &\quad + \gamma^{t+1}\sum_{s'\in\setS} \nu_{t+1,s'}  \left(P(s'|s, a, x_t(\mu,\pi)) \mu_t(s) + \mathbb{E}\left[\frac{\partial P(s'|\cdot, x_t(\mu,\pi))}{\partial x_t(s,a)}\right]\mu_t(s)\right), \nonumber\\
        0 &= \gamma^t \sum_{a\in\setA} \pi_t(s,a) R(s,a,x_t(\mu,\pi))  - \gamma^t \nu_{t,s} \\
        &\quad + \gamma^{t+1}\sum_{s'\in\setS} \nu_{t+1,s'} \sum_{a\in\setA} \left(P(s'|s, a, x_t(\mu,\pi)) \pi_t(s,a) + \mathbb{E}\left[\frac{\partial P(s'|\cdot,x_t(\mu,\pi))}{\partial x_t(s,a)}\right]\pi_t(s,a)\right), \nonumber
    \end{align}
\end{subequations}
Note that the partial derivatives $\partial \Phi(x(\mu, \pi))/\partial \pi_t(s,a)$ and $\partial \Phi(x(\mu, \pi))/\partial \mu_t(s)$ are now decoupled from future values thanks to the constraints and multipliers. Below, we will show the multipliers are closely related to the state and state-action values at equilibrium. 

Since $\gamma>0$ and due to Assumption~\ref{ass:transition}, when $\mu_t(s) > 0$, the above equations are further reduced to
\begin{align}
    \lambda_{t,s}  & = R(s,a,x_t(\mu,\pi)) + \eta_{t,sa} + \gamma \sum_{s'\in\setS} P(s'|s, a, x_t(\mu,\pi)) \nu_{t+1,s'} \label{eq:KKT_lambda_ass2}
    \\
    \nu_{t,s}  & = \sum_{a\in\setA}\pi_t(s,a) \left[R(s,a,x_t(\mu,\pi)) +  \gamma \sum_{s'\in\setS} P(s'|s,a,x_t(\mu,\pi))  \nu_{t+1, s'}\right] \label{eq:KKT_nu_ass2}.
\end{align}
On the other hand, KKT conditions~\eqref{eq:simplified_KKT_lambda} are clearly satisfied when $\mu_t(s) = 0$. Thus, we focus on the case when $\mu_t(s) > 0$.

To formally establish the equivalence between multipliers and value functions, we employ 
backward induction to prove that 
\begin{align}
    \lambda_{t,s} = Q^*_t(s,x(\mu,\pi)),
    \quad
    \nu_{t,s} = V_t(s,x(\mu,\pi)),
\end{align}
which yields the optimality condition \eqref{eq:def_equ_Q}. 

Without violating any constraints, we can set $\nu_{T,s}=V_T(s,x(\mu,\pi))=0, \forall s\in \setS$. Then, it is easy to verify that for time $T-1$, we have
\begin{align}
    \lambda_{T-1,s}  & = R(s,a,x_{T-1}(\mu,\pi)) + \eta_{T-1,sa} = Q_{T-1}(s,a,x(\mu,\pi)) + \eta_{T-1,sa},\\
    \nu_{T-1,s} &= \sum_{a\in\setA}\pi_{T-1}(s,a) R(s,a,x_{T-1}(\mu,\pi))  = V_{T-1}(s,x(\mu,\pi)).
\end{align}
Combining with the complementary slackness \eqref{eq:KKT_comple} and \eqref{eq:KKT_nonneg}, we have 
\begin{align}
    \lambda_{T-1,s} &= Q_{T-1}(s,a,x(\mu,\pi)), & \text{ if } \pi_{T-1}(s,a) > 0,\\
    \lambda_{T-1,s} &\geq Q_{T-1}(s,a,x(\mu,\pi)), & \text{ if } \pi_{T-1}(s,a) = 0,
\end{align}
which implies $\lambda_{T-1,s} = Q^*_{T-1}(s,x(\mu,\pi))$.

Next, assuming the result holds for some $t+1\leq T$, we have for time $t$
\begin{align}\label{eq:induction_eq}
    \lambda_{t,s}  & = R(s,a,x_t(\mu,\pi)) + \eta_{t,sa} + \gamma \sum_{s'\in\setS} P(s'|s, a, x_t(\mu,\pi))V_{t+1}(s',x(\mu,\pi)) = Q_t(s,a,x(\mu,\pi)) + \eta_{t,sa},\\
    \nu_{t,s} &= \sum_{a\in\setA}\pi_t(s,a) \left[R(s,a,x_t(\mu,\pi)) +  \gamma \sum_{s'\in\setS} P(s'|s,a,x_t(\mu,\pi)) V_{t+1}(s',x(\mu,\pi))\right] = V_t(s, x(\mu, \pi)).
\end{align}
Same as the proof for $T-1$, the complementary slackness \eqref{eq:KKT_comple} and \eqref{eq:KKT_nonneg} imply $\lambda_{t,s} = Q^*_t(s,x(\mu,\pi))$.
It thus follows that $\lambda_{t,s} = Q^*_t(s,x(\mu,\pi))$ and the optimality condition \eqref{eq:def_equ_Q} holds for $t$.

In sum, a tuple $(\mu, \pi)$ satisfying the KKT conditions \eqref{eq:KKT}  also satisfies the equilibrium conditions \eqref{eq:def_equ}. Thus, $(\mu, \pi)$ also describes an equilibrium of MPPG with potential $\Phi(x(\mu, \pi))$.

\section{Proof of Proposition~\ref{prop:equ_identical_Q}}\label{app:equ_identical_Q}

The first statement is easy to prove because any KKT point satisfies Eq.~\eqref{eq:eps_equ_mu} and Eq.~\eqref{eq:eps_equ_BR} also holds due to the bounded reward and finite horizon. 

To prove the second statement, we first rewrite the KKT condition derived in Appendix~\ref{app:prop_optm_MPPG} as
\begin{align}
    \lambda_{t,s}  & = R(s,a,\tilde x_t) + \gamma \sum_{s'\in\setS}
   \left(P(s'|s,a,\tilde x_t)
   + \mathbb{E}\!\left[\frac{\partial P(s'|\cdot,\tilde x_t)}
     {\partial x_t(s,a)}\right]\right) \nu_{t+1,s'}, && \tilde \pi_t(s,a) > 0, \label{eq:equ_identical_KKT_1}\\
     \lambda_{t,s}  &\geq R(s,a,\tilde x_t) + \gamma \sum_{s'\in\setS}
   \left(P(s'|s,a,\tilde x_t)
   + \mathbb{E}\!\left[\frac{\partial P(s'|\cdot,\tilde x_t)}
     {\partial \tilde x_t(s,a)}\right]\right) \nu_{t+1,s'}, && \tilde \pi_t(s,a) = 0.\label{eq:equ_identical_KKT_2}
\end{align}

Suppose the equality $\nu_{t+1,s'}= V_{t+1}(s',\tilde x) + \hat \epsilon_{t+1}$ for some constant $\hat\epsilon_{t+1}$ holds. Then, Eq.~\eqref{eq:equ_identical_KKT_1} can be rewritten as 
\begin{align}
    \lambda_{t,s}  & = R(s,a,\tilde x_t) + \gamma \sum_{s'\in\setS}
   P(s'|s,a,\tilde x_t)( V_{t+1}(s',\tilde x) + \hat \epsilon_{t+1}) +  \gamma \sum_{s'\in\setS}
   \mathbb{E}\!\left[\frac{\partial P(s'|\cdot,\tilde x_t)}
     {\partial x_t(s,a)}\right] \nu_{t+1,s'}\nonumber\\
     &= R(s,a,\tilde x_t) + \gamma \sum_{s'\in\setS}
   P(s'|s,a,\tilde x_t) V_{t+1}(s',\tilde x) + \gamma (\hat \epsilon_{t+1} + \bar \epsilon_{t}),
\end{align}
where the second equality is due to the condition~\eqref{eq:local_stationarity}. 

Applying the same derivation to Eq.~\eqref{eq:equ_identical_KKT_2} leads to
\begin{align}
    \lambda_{t,s}  & = R(s,a,\tilde x_t) + \gamma \sum_{s'\in\setS}
   P(s'|s,a,\tilde x_t) V_{t+1}(s',\tilde x) + \gamma (\hat \epsilon_{t+1} + \bar \epsilon_{t}), && \tilde \pi_t(s,a) > 0,\\
     \lambda_{t,s}  &\geq R(s,a,\tilde x_t) + \gamma \sum_{s'\in\setS}
   P(s'|s,a,\tilde x_t) V_{t+1}(s',\tilde x) + \gamma (\hat \epsilon_{t+1} + \bar \epsilon_{t}), && \tilde \pi_t(s,a) = 0,
\end{align}
which is equivalent to the optimality condition \eqref{eq:def_equ_pi_BR} of equilibrium because the constant shift $\gamma (\hat \epsilon_{t+1} + \bar \epsilon_{t})$ does not change the optimality of $\tilde \pi$. 

Below, we prove the result $\nu_{t+1,s'}= V_{t+1}(s',\tilde x) + \hat \epsilon_{t+1}$ by backward induction. At the terminal step $t=T$, the result naturally holds with $\hat \epsilon_{T} = 0,\forall s$ as $\nu_{T,s} = V_{T}(s,\tilde x)=0,\;\forall s$. Suppose the result holds for $t+1$, the other KKT condition derived in Appendix~\ref{app:prop_optm_MPPG} reads 
\begin{align}
    \nu_{t,s}  &= \sum_{a\in\setA}\tilde\pi_t(s,a) \left[R(s,a,\tilde x_t) +  \gamma \sum_{s'\in\setS}
   \left(P(s'|s,a,\tilde x_t)
   + \mathbb{E}\!\left[\frac{\partial P(s'|\cdot,\tilde x_t)}
     {\partial \tilde x_t(s,a)}\right]\right) \nu_{t+1, s'}\right]\nonumber\\
     &= \sum_{a\in\setA}\tilde\pi_t(s,a) \left[R(s,a,\tilde x_t) +  \gamma \sum_{s'\in\setS}
   P(s'|s,a,\tilde x_t)(V_{t+1}(s',\tilde x) + \hat \epsilon_{t+1})\right] + \gamma \sum_{a\in\setA}\tilde\pi_t(s,a) \sum_{s'\in\setS} \mathbb{E}\!\left[\frac{\partial P(s'|\cdot,\tilde x_t)}{\partial \tilde x_t(s,a)}\right] \nu_{t+1, s'}\nonumber\\
     &= V_t(s,\tilde x) + \gamma(\hat \epsilon_{t+1} + \bar \epsilon_{t}). 
\end{align}
The last equality is again due to condition~\eqref{eq:local_stationarity}. Hence, the result holds with $\hat \epsilon_{t} = \gamma(\hat \epsilon_{t+1} + \bar \epsilon_{t})$.

\section{Proof of Proposition~\ref{theo:priori bound}}
\label{app:theo_priori_bound}

The KKT conditions derived in Appendix~\ref{app:prop_optm_MPPG} do not rely on Assumption~\ref{ass:transition} and remain necessary. Thus, Eq.~\eqref{eq:KKT_eq_mu} still corresponds to the equilibrium condition about state distribution dynamics, or equivalently, Eq.~\eqref{eq:eps_equ_mu}, and we need to prove that Condition~\eqref{eq:eps_equ_BR} holds with a priori bound defined in Eq.~\eqref{eq:epsilon_apriori_cases}. 
Throughout the proof, we simplify the notation by omitting $(\mu,\pi)$ in reward, transition, and value functions, e.g., $R_t(s,a)=R(s,a,x_t)$.
We also denote the expected marginal transition effect by
$D_t(s'\mid s,a)
=\mathbb{E}\left[
\frac{\partial P(s'\mid\cdot,x_t)}
{\partial x_t(s,a)}
\right]$.

According to the KKT conditions shown in Appendix~\ref{app:prop_optm_MPPG}, at state $s$ with $\mu_t(s)>0$, we have
\begin{align}
    \lambda_{t,s}  & = R_t(s,a) + \gamma \sum_{s'\in\setS}
   \left(P_t(s'|s,a)
   + D_t(s'|s,a)\right) \nu_{t+1,s'} ,\; \forall a\in\setA_+(t,s;\pi)\label{eq:KKT_lambda_w/_partial}\\
    \nu_{t,s}  & = \sum_{a\in\setA}\pi_t(s,a) \left[R_t(s,a) +  \gamma \sum_{s'\in\setS}
   \left(P_t(s'|s,a)
   + D_t(s'|s,a)\right) \nu_{t+1, s'}\right] \label{eq:KKT_nu_w/_partial}
\end{align}
The dual variables $\lambda_{t,s}$ and $\nu_{t,s}$ correspond  to pseudo state-action value $\tilde{Q}_t(s,a)$ and pseudo value function $\tilde{V}_t(s)$, respectively. Accordingly, for any two actions $a_1,a_2\in\setA_+(t,s;\pi)$, we have 
$\tilde Q_t(s,a_1)=\tilde Q_t(s,a_2)$.

Write $P_t$ and $D_t$ as matrices with rows indexed
by $(s,a)$ and columns indexed by $s'$. 
Then, Eqs.~\eqref{eq:KKT_lambda_w/_partial} and \eqref{eq:KKT_nu_w/_partial} can be rewritten in matrix form:
\begin{subequations}
    \begin{align}
\tilde Q_t
&=
R_t+\gamma(P_t+D_t)\tilde V_{t+1},
\label{eq:apriori_pseudo_Q}\\
\tilde V_t(s)
&=
\sum_{a\in\setA}\pi_t(s,a)\tilde Q_t(s,a),
\quad
\tilde V_T(\cdot)=0.
\label{eq:apriori_pseudo_V}
\end{align}
\end{subequations}
On the other hand, under the same $(\mu, \pi)$ the true value functions satisfy
\begin{subequations}
    \begin{align}
Q_t
&=
R_t+\gamma P_tV_{t+1},
\label{eq:apriori_true_Q}\\
V_t(s)
&=
\sum_{a\in\setA}\pi_t(s,a)Q_t(s,a),
\quad V_T(\cdot)=0.
\label{eq:apriori_true_V}
\end{align}
\end{subequations}

Define the Q-value residual as 
$\delta_t(s,a)=\tilde Q_t(s,a)-Q_t(s,a)$. Since any two active actions $a_1,a_2$ satisfy $\widetilde Q_t(s,a_1)=\widetilde Q_t(s,a_2)$, we have
\begin{align}
|Q_t(s,a_1)-Q_t(s,a_2)|=
|\delta_t(s,a_1)-\delta_t(s,a_2)|\leq\min\left\{
\operatorname{span}(Q_t),
\operatorname{span}(\delta_t)
\right\},
\label{eq:apriori_gap_reduction}
\end{align}
where $\text{span}(\cdot)$ denotes the span of a finite-dimensional vector, i.e., $\operatorname{span}(e)=\max_i e_i-\min_i e_i$.

We first derive uniform bounds for $\operatorname{span}(Q_t)$.
For any $s\in\setS, a\in\setA$, 
\begin{subequations}
    \begin{align}
\sum_{s'\in\setS}D_t(s'\mid s,a)= \sum_{\bar s \in\setS, \bar a\in\setA}x_t(\bar s,\bar a) \frac{\partial 
\sum_{s'\in\setS} P_t(s'\mid\bar s,\bar a)
}{\partial x_t(s,a)}
=0,
\label{eq:apriori_D_zero_sum}
\end{align}
\begin{align}
\sum_{s'\in\setS}|D_t(s'\mid s,a)|\leq
\sum_{\bar s\in\setS, \bar a\in\setA}
x_t(\bar s,\bar a)
\sum_{s'\in\setS}
\left|
\frac{\partial P_t(s'\mid\bar s,\bar a)}
{\partial x_t(s,a)}
\right|\leq |\setS|L_P.
\label{eq:apriori_D_row_bound}
\end{align}
\end{subequations}

For any state value vector $v$, define the mean value as
$\bar v=\frac12
\left(\max_{s}v(s)+\min_{s}v(s)\right)$, then
$|v(s')-\bar v|
\leq\frac12\operatorname{span}(v),\forall s'\in\setS$.
Due to Eq.~\eqref{eq:apriori_D_zero_sum}, we have $(D_t\bar{v})(s,a) = \sum_{s'\in\setS} D_t(s'\mid s,a)\bar v = 0$. Hence, it yields
\begin{align}
|(D_tv)(s,a)| = |(D_tv)(s,a) - (D_t\bar{v})(s,a) | 
&=
\left|
\sum_{s'\in\setS}
D_t(s'\mid s,a)\bigl(v(s')-\bar v\bigr)
\right|
\nonumber\\
&\leq
\frac12\operatorname{span}(v)
\sum_{s'\in\setS}|D_t(s'\mid s,a)|
\leq
\frac{|\setS|L_P}{2}\operatorname{span}(v),
\end{align}
which leads to $\operatorname{span}(D_tv)\leq|\setS|L_P\operatorname{span}(v)$.
Similarly, we have
\begin{align}
\operatorname{span}(P_tv)
&=
\max_{\substack{s_1,s_2\in\setS\\a_1,a_2\in\setA}}
\left|
\sum_{s'\in\setS}
\bigl[
P_t(s'\mid s_1,a_1)-P_t(s'\mid s_2,a_2)
\bigr]
\bigl(v(s')-\bar v\bigr)
\right|
\nonumber\\
&\leq
\frac12\operatorname{span}(v)
\max_{\substack{s_1,s_2\in\setS\\a_1,a_2\in\setA}}
\sum_{s'\in\setS}
\left|
P_t(s'\mid s_1,a_1)-P_t(s'\mid s_2,a_2)
\right|\leq
\theta\operatorname{span}(v).
\label{eq:apriori_P_span_bound}
\end{align}
Besides, Eqs.~\eqref{eq:stage reward} and \eqref{eq:def_V} imply $\operatorname{span}(R_t) \leq\Delta r$ and $\operatorname{span}(V_t) \leq\operatorname{span}(Q_t)$. 
Combining these bounds with Eq.~\eqref{eq:apriori_true_Q},  
we obtain
\begin{align}
\operatorname{span}(Q_t)
\leq
\Delta r+\gamma\theta
\operatorname{span}(V_{t+1})\leq
\Delta r+\gamma\theta
\operatorname{span}(Q_{t+1}),
\end{align}
At the final stage,
$Q_{T-1}=R_{T-1}$, and hence
$\operatorname{span}(Q_{T-1})\leq\Delta r$.
Backward induction gives
\begin{equation}
\operatorname{span}(Q_t)
\leq
\Delta r
\sum_{k=0}^{T-t-1}(\gamma\theta)^k\leq
\Delta r\sum_{k=0}^{T-1}(\gamma\theta)^k,
\quad t\in\setT,
\label{eq:apriori_true_value_stage_bound}
\end{equation}
where the last inequality holds because all summands are nonnegative.

Next, we derive the bound for $\operatorname{span}(\delta_t)$. First, rewrite 
$\delta_t$ as
\begin{align}
\delta_t
=\gamma D_tV_{t+1}
+\gamma P_t(\tilde V_{t+1}-V_{t+1})+
\gamma D_t(\tilde V_{t+1}-V_{t+1}).
\label{eq:apriori_delta_recursion}
\end{align}
Then, applying the bounds of $\operatorname{span}(D_tv)$ and $\operatorname{span}(P_tv)$ yields
\begin{align}
\operatorname{span}(\delta_t)
&\leq
\gamma|\setS|L_P
\operatorname{span}(V_{t+1})+
\gamma\bigl(\theta+|\setS|L_P\bigr)
\operatorname{span}(\tilde V_{t+1}-V_{t+1})
\nonumber\\
&\leq
\gamma|\setS|L_P
\operatorname{span}(V_{t+1})
+
\gamma\bigl(\theta+|\setS|L_P\bigr)
\operatorname{span}(\delta_{t+1}),
\label{eq:apriori_delta_span_recursion}
\end{align}
To see the last inequality, recall that equations~\eqref{eq:apriori_pseudo_V}
and~\eqref{eq:apriori_true_V} mean $\tilde V_t(s)-V_t(s)=
\sum_{a\in\setA}\pi_t(s,a)\delta_t(s,a)$,
implying $\operatorname{span}(\tilde V_t-V_t)
\leq\operatorname{span}(\delta_t)$.

Since $\tilde Q_{T-1}=Q_{T-1}=R_{T-1}$, we have $\operatorname{span}(\delta_{T-1})=0$.
Thus, backward induction gives
\begin{align}
\operatorname{span}(\delta_t)
&\leq
\gamma|\setS|L_P
\sum_{j=0}^{T-t-2}
\bigl[\gamma(\theta+|\setS|L_P)\bigr]^j
\operatorname{span}(V_{t+j+1})
\nonumber\\
&\leq
\Delta r \gamma|\setS|L_P
\sum_{j=0}^{T-t-2}
\bigl[\gamma(\theta+|\setS|L_P)\bigr]^j
\sum_{k=0}^{T-t-j-2}(\gamma\theta)^k,
\label{eq:apriori_error_double_sum}
\end{align}
where the second inequality is due to $\operatorname{span}(V_t) \leq\operatorname{span}(Q_t)$ and Eq.~\eqref{eq:apriori_true_value_stage_bound}.

Set $m=j+k+1$. The domain
$j,k\geq0$, $j+k\leq T-t-2$ then becomes
$1\leq m\leq T-t-1$, $0\leq j\leq m-1$, with
$k=m-1-j$. Accordingly, the bound \eqref{eq:apriori_error_double_sum} is further reduced to
\begin{align}
\operatorname{span}(\delta_t)
&\leq \Delta r\gamma|\setS|L_P
\sum_{j=0}^{T-t-2}
\bigl[\gamma(\theta+|\setS|L_P)\bigr]^j
\sum_{k=0}^{T-t-j-2}(\gamma\theta)^k
= \Delta r|\setS|L_P
\sum_{j=0}^{T-t-2}
\sum_{k=0}^{T-t-j-2}
\gamma^{j+k+1}
(\theta+|\setS|L_P)^j\theta^k
\nonumber\\
&= \Delta r
|\setS|L_P
\sum_{m=1}^{T-t-1}
\sum_{j=0}^{m-1}
\gamma^m
(\theta+|\setS|L_P)^j\theta^{m-1-j}
= \Delta r\sum_{m=1}^{T-t-1}\gamma^m
\left[
|\setS|L_P
\sum_{j=0}^{m-1}
(\theta+|\setS|L_P)^j\theta^{m-1-j}
\right]
\nonumber\\
&= \Delta r \sum_{m=0}^{T-t-1}\gamma^m
\left[
(\theta+|\setS|L_P)^m-\theta^m
\right].
\label{eq:apriori_error_stage_bound}
\end{align}
where the last equality follows from the property
$a^m-b^m=(a-b)\sum_{j=0}^{m-1}a^jb^{m-1-j}$.
Since all summands are nonnegative, it is safe to drop the index $t$ and derive
\begin{equation}
\operatorname{span}(\delta_t)
\leq
\Delta r
\sum_{m=0}^{T-1}
\gamma^m
\left[
(\theta+|\setS|L_P)^m-\theta^m
\right],
\qquad t\in\setT.
\label{eq:apriori_error_uniform_bound}
\end{equation}
Finally, combining
Eqs.~\eqref{eq:apriori_gap_reduction},
\eqref{eq:apriori_true_value_stage_bound},
and~\eqref{eq:apriori_error_uniform_bound},
we obtain
\begin{align}
|Q_t(s,a_1)-Q_t(s,a_2)|
\leq
\Delta r
\min\Bigg\{
\sum_{k=0}^{T-1}(\gamma\theta)^k,
\sum_{k=0}^{T-1}\gamma^k
\left[
(\theta+|\setS|L_P)^k-\theta^k
\right]
\Bigg\}
\label{eq:apriori_finite_sum_bound}
\end{align}
for every pair of active
actions $a_1,a_2$.
This concludes the proof.

\section{Proof of Proposition~\ref{theo:posteriori_bound}}
\label{app:posteriori_bound}

Fix the KKT point $(\mu,\pi)$ and its associated multipliers.
Throughout the proof, we employ the simplified notations in Appendix~\ref{app:theo_priori_bound}. Besides, we focus on states $s$ at time $t$ with $\mu_t(s)>0$, and all action
comparisons are restricted to $\mathcal A_+(t,s;\pi)$.
The KKT conditions in Appendix~\ref{app:prop_optm_MPPG} imply
\begin{equation}
\label{eq:post_dual_identity}
\lambda_{t,s}=\sum_{a\in\setA_+(t,s;\pi)}
\pi_t(s,a)\lambda_{t,s}=\nu_{t,s}.
\end{equation}
Combining it with the terminal condition $\nu_{T,\cdot}=V_T(\cdot,\mu,\pi)=0$ yields
\begin{equation}
\label{eq:posteriori_terminal}
\mathcal{B}_{T-1}(s,a)=0,
\quad
\nu_{T-1,s}=V_{T-1}(s).
\end{equation}
Averaging Eq.~\eqref{eq:posteriori_residual} with respect to
$\pi_t(s,\cdot)$ and applying
Eq.~\eqref{eq:post_dual_identity} gives
\begin{align}
\nu_{t,s}
=
\overline{\mathcal{B}}_t(s)
+
\sum_{a\in\mathcal A_+(t,s;\pi)}
\pi_t(s,a)
\left[
R_t(s,a)
+
\gamma\sum_{s'\in\setS}
P_t(s'\mid s,a)\nu_{t+1,s'}
\right].
\end{align}
Define the value error
as $e_t(s)=\nu_{t,s}-V_t(s,\mu,\pi)$, which yields
\begin{equation}
\label{eq:posteriori_value_error}
e_t(s)
=
\overline{\mathcal{B}}_t(s)
+
\gamma
\sum_{a\in\mathcal A_+(t,s;\pi)}
\pi_t(s,a)
\sum_{s'\in\setS}
P_t(s'\mid s,a)e_{t+1}(s'), \; \forall t\in\setT.
\end{equation}
For the policy-averaged transition, the pairwise total variations are also
bounded by the same uniform Dobrushin coefficient $\theta$.
Thus,
\begin{align}
\operatorname{span}\left(
\sum_{a\in\mathcal A_+(t,s;\pi)}
\pi_t(s,a)
\sum_{s'\in\setS}P_t(s'\mid s,a)e_{t+1}(s')
\right)
\le
\theta\operatorname{span}(e_{t+1}),\quad 0\le t\le T-2.
\end{align}
Combining this inequality with
Eq.~\eqref{eq:posteriori_value_error} yields
\begin{align}
\operatorname{span}(e_t)
\le
\operatorname{span}(\overline{\mathcal{B}}_t)
+
\gamma\theta\operatorname{span}(e_{t+1})
\le
\mathcal{B}_{\mathrm{avg}}
+
\gamma\theta\operatorname{span}(e_{t+1}),
\quad 0\le t\le T-2.
\label{eq:posteriori_span_recursion}
\end{align}
Since $e_{T-1}=0$, backward induction gives
\begin{equation}
\label{eq:posteriori_value_error_bound}
\operatorname{span}(e_t)
\le
\mathcal{B}_{\mathrm{avg}}
\sum_{j=0}^{T-t-2}(\gamma\theta)^j,
\; t\in\setT.
\end{equation}

Next, on an active action $a$, the true state-action value function satisfies
\begin{align}
Q_t(s,a)
&=
R_t(s,a)+
\gamma\sum_{s'\in\setS}
P_t(s'\mid s,a)V_{t+1}(s')
=
\lambda_{t,s}-\mathcal{B}_t(s,a)-
\gamma\sum_{s'\in\setS}
P_t(s'\mid s,a)e_{t+1}(s').
\label{eq:posteriori_Q_representation}
\end{align}
For two active actions at the same state, it's derived that
\begin{align}
\left|Q_t(s,a_1)-Q_t(s,a_2)\right|
&\le
\left|\mathcal{B}_t(s,a_1)-\mathcal{B}_t(s,a_2)\right|+
\gamma\left|\sum_{s'\in\setS}\left[P_t(s'\mid s,a_1)-P_t(s'\mid s,a_2)\right]e_{t+1}(s')\right|
\nonumber\\
&\le\mathcal{B}_{\mathrm{act}}+\gamma\theta\operatorname{span}(e_{t+1})\le
\mathcal{B}_{\mathrm{act}}+\mathcal{B}_{\mathrm{avg}}
\sum_{k=1}^{T-t-2}(\gamma\theta)^k\nonumber\\
&\le\mathcal{B}_{\mathrm{act}}+\mathcal{B}_{\mathrm{avg}}
\sum_{k=1}^{T-2}(\gamma\theta)^k.
\label{eq:posteriori_uniform_gap}
\end{align}
The last inequality follows from the nonnegativity of each
summand.
Taking the minimum with the a priori bound completes the proof.

\section{Lemma~\ref{lemma:unique_V} and proof}\label{app:Lemma_unique_V}

\begin{lemma}\label{lemma:unique_V}
    The KKT conditions Eqs.~\eqref{eq:KKT_grad_mu_sub} and \eqref{eq:KKT_eq_mu_sub} admit a unique solution of $\mu,\nu$.
\end{lemma}

\begin{proof}
    First note that the state distribution $\mu$ is uniquely determined by Eq.~\eqref{eq:KKT_eq_mu_sub} given policy $\pi$ and initial state distribution $\mu_0$. Accordingly, we define constant vector and matrices
    \begin{align}
    \bar R_t(s)&:=\sum_{a\in A}\pi_t(s,a)\,R(s,a,x_t),\\
    \bar P_t(s,s')&:=\sum_{a\in A}\pi_t(s,a)\left( P(s'|s,a,x_t)+ \mathbb{E}\!\left[\frac{\partial P(s'|\cdot,x_t)}{\partial x_t(s,a)}\right]\right).
    \end{align}
    Then Eq.~\eqref{eq:KKT_grad_mu_sub} is reduced to $\nu_t=\bar R_t+\gamma\,\bar P_t\,\nu_{t+1}$ for any $t\in\setT$.
    By setting the terminal value $\nu_T=0$, we obtain $\nu$ uniquely by backward induction. 
\end{proof}

\section{Proof of Theorem~\ref{prop:convergence}}\label{app:prop_convergence}
 
We first show that every limit point of the policy
sequence $\{\pi^{(i)}\}$ generated by the PG algorithm satisfies the KKT
conditions of the master problem \eqref{eq:optm_mas}.
For each time and state, the projection onto a convex set satisfies
\begin{align}\label{eq:prop_convergence_proj}
\big\langle
&\pi_t^{(i)}(s,\cdot)
+\eta_{t,s}^{(i)}
\nabla F(\pi^{(i)})_{t,s,\cdot}
-\pi_t^{(i+1)}(s,\cdot),\;\pi_t^{(i)}(s,\cdot)-\pi_t^{(i+1)}(s,\cdot)
\big\rangle
\leq0.
\end{align}
Let $d^{(i)}_{t,s}=\pi_t^{(i+1)}(s,\cdot)-\pi_t^{(i)}(s,\cdot)$. Then Eq.~\eqref{eq:prop_convergence_proj} is rewritten as
\begin{align}
    & \langle \eta_{t,s}^{(i)}
\nabla F(\pi^{(i)})_{t,s,\cdot} - d^{(i)}_{t,s}, \; - d^{(i)}_{t,s}\rangle \leq 0,\\
\Rightarrow \quad & \langle \nabla F(\pi^{(i)})_{t,s,\cdot},  d^{(i)}_{t,s}\rangle \geq \frac{1}{\eta_{t,s}^{(i)}} \|d_t^{(i)}(s,\cdot)\|_2^2.
\end{align}
Summing over all times and states gives 
\begin{align}
    \sum_{t,s}\langle \nabla F(\pi^{(i)})_{t,s,\cdot},  d^{(i)}_{t,s}\rangle = \langle \nabla F(\pi^{(i)}),  d^{(i)}\rangle \geq \sum_{t,s} \frac{1}{\eta_{t,s}^{(i)}} \|d_t^{(i)}(s,\cdot)\|_2^2 \geq \frac{1}{\bar\beta}\|d^{(i)}\|_2^2,
\end{align}
where the last inequality is due to $0<\eta_{t,s}^{(i)}\leq\bar\beta$.

The Lipschitz continuity of the gradient therefore implies
\begin{align}
F(\pi^{(i+1)})-F(\pi^{(i)})
\geq
\left\langle\nabla F(\pi^{(i)}),d^{(i)}\right\rangle
-\frac{L_F}{2}\|d^{(i)}\|_2^2\geq
\left(\frac{1}{\bar\beta}-\frac{L_F}{2}\right)
\|d^{(i)}\|_2^2 > 0,
\end{align}
where the last inequality holds because $\bar\beta<\frac{2}{L_F}$. 

Since the feasible policy set $\Pi$ is compact and $F$ is continuous, $F$ is bounded above on $\Pi$. Therefore, the monotonically increasing sequence
$\{F(\pi^i)\}$ converges to a finite limit $F_\infty$. 
Summing the improvement over $i=0,\ldots,K-1$ gives
\begin{align}
\left(\frac{1}{\bar\beta}-\frac{L_F}{2}\right)
\sum_{i=0}^{K-1}\|\pi^{i+1}-\pi^i\|_2^2
\leq F(\pi^K)-F(\pi^0)
\leq F_\infty-F(\pi^0).
\end{align}
Since $\bar \beta < \frac{2}{L_F}$, letting $K\to\infty$ yields $\sum_{i=0}^{\infty}\|\pi^{i+1}-\pi^i\|_2^2<\infty$,
which dictates that $\lim_{i \to \infty} |\pi^{(i+1)} - \pi^{(i)}|= 0$. As $\Pi$ is compact, ${\pi^{(i)}}$ has at least one limit point $\pi^*$. 
Let $\pi^*$ be an arbitrary limit point and
$\{\pi^{(i_k)}\}$ be a subsequence such that $\pi^{(i_k)}\to\pi^*$.
It implies that $\pi^{(i_k+1)}\to\pi^*$ and 
\begin{align}
\pi_t^*(s,\cdot)
=
\operatorname{Proj}_{\Delta(\mathcal A)}
\left[
\pi_t^*(s,\cdot)
+
\eta_{t,s}
\nabla F(\pi^*)_{t,s,\cdot}
\right].
\end{align}
Applying the projection inequality again, we can derive
\begin{align}
\left\langle
\nabla F(\pi^*),\pi-\pi^*
\right\rangle
\leq0,
\quad \forall \pi\in\Pi,
\end{align}
which corresponds to the first-order stationarity condition of the master problem~\eqref{eq:optm_mas}. Since $\Pi$ is a convex polytope, $\pi^*$ also satisfies the KKT conditions of Problem~\eqref{eq:optm_mas}~\citep{facchinei2003finite}. Accordingly, $\nu^{\pi^*}$ and $\mu^{\pi^*}$ satisfy
Eqs.~\eqref{eq:KKT_grad_mu_sub} and \eqref{eq:KKT_eq_mu_sub}, respectively. Then all the KKT conditions of Problem~\eqref{eq:optm_MPPG} are satisfied, and we complete the proof.

\section{Proof of Proposition~\ref{prop:PUME}}\label{app:PUME}

Consider a finite-horizon version of PUME with exogenous transition and separable link utilities, which can be embedded in the original infinite-horizon formulation by defining the absorbing state with the terminal time. For a feasible profile $(\mu,\pi)$ with induced flow $x$, define
\begin{align}
R(s,a,x_t)
&=
r(s,a,x_t)
-\frac{\partial H_s^*(\pi_t(s,\cdot))}
       {\partial\pi_t(s,a)}
+c_t(s,x_t)\label{eq:PUME_reward},\\
c_t(s,x_t)
&=
\pi_t(s,\cdot)^\top\nabla H_s^*(\pi_t(s,\cdot))
-H_s^*(\pi_t(s,\cdot)). \label{eq:PUME_correction}
\end{align}
The corresponding value function $V$ and $Q$-value function satisfy
\begin{align}
Q_t(s,a,x)
&=
R(s,a,x_t)
+\gamma\sum_{s'\in\setS}P(s'\mid s,a)V_{t+1}(s',x),\label{eq:PUME_Q}\\
V_t(s,x)
&=
\sum_{a\in\setA}\pi_t(s,a)Q_t(s,a,x)
\nonumber\\
&=
\sum_{a\in\setA}\pi_t(s,a)
\left[
r(s,a,x_t)
+\gamma\sum_{s'\in\setS}P(s'\mid s,a)V_{t+1}(s',x)
\right]
-H_s^*(\pi_t(s,\cdot)),
\label{eq:PUME_value}
\end{align}
with $V_T(s,x)=0, \forall s\in\setS$.

Let $(\mu^*,\pi^*)$ be an equilibrium of the finite-horizon
PUME model, with induced flow $x^*$.
At each time and state, $\pi_t^*(s,\cdot)$ solves
\begin{subequations}\label{eq:pume}
\begin{align}
\max_{\pi_t(s,\cdot)}\quad
&
\sum_{a\in\setA}\pi_t(s,a)
\left[
r(s,a,x_t^*)
+\gamma\sum_{s'\in\setS}
P(s'\mid s,a)V_{t+1}(s',x^*)
\right]
-H_s^*(\pi_t(s,\cdot))
\\
\text{s.t.}\quad
&\sum_{a\in\setA}\pi_t(s,a)=1,
\label{eq:pume_normalization}\\
&\pi_t(s,a)\geq0,\qquad \forall a\in\setA.
\label{eq:pume_nonnegativity}
\end{align}
\end{subequations}
Given $x^*$ and convex $H_s^*$, it is easy to show Problem~\eqref{eq:pume} is convex and thus the KKT conditions are also sufficient:
\begin{align}
\lambda_{t,s}
&=
r(s,a,x_t^*)
-\frac{\partial H_s^*(\pi_t^*(s,\cdot))}
       {\partial\pi_t(s,a)}
+\gamma\sum_{s'\in\setS}P(s'\mid s,a)V_{t+1}(s',x^*)
+\eta_{t,sa}, \label{eq:PUME_stationarity}\\
\eta_{t,sa}&\geq0,
\qquad
\eta_{t,sa}\pi_t^*(s,a)=0,
\label{eq:PUME_slackness}
\end{align}
where $\lambda_{t,s}$ and $\eta_{t,sa}$ are the
multipliers for  constraints~\eqref{eq:pume_normalization} and ~\eqref{eq:pume_nonnegativity}, respectively.

With $c_t$ and $Q_t$ defined in Eqs.~\eqref{eq:PUME_correction} and ~\eqref{eq:PUME_Q}, Eq.~\eqref{eq:PUME_stationarity} can be written as
\begin{align}\label{eq:PUME_maxQ}
    \lambda_{t,s}
&=
Q_t(s,a,x^*)
-c_t(s,x^*_t)
+\eta_{t,sa}.
\end{align}
Since $c_t(s,x^*_t)$ is action-independent, combining Eq.~\eqref{eq:PUME_maxQ} and complementary slackness~\eqref{eq:PUME_slackness} therefore implies that every active
action $a$ satisfies: 
\begin{equation}
\pi_t^*(s,a)
\left[Q_t(s,a,x^*)-Q_t^*(s,x^*)\right]=0,
\end{equation}
which gives
the equilibrium conditions in Eq.~\eqref{eq:def_equ_Q}. Since $(\mu^*,\pi^*)$ also satisfies the population dynamics
in Eq.~\eqref{eq:def_equ_mu}, we prove $(\mu^*,\pi^*)$ corresponds to an equilibrium of MPPG.

For $\mu_t(s)>0$, using
$\mu_t(s)=\sum_b x_t(s,b)$ and
$\pi_t(s,b)=x_t(s,b)/\mu_t(s)$ gives
\begin{align}
\frac{\partial \pi_t(s,b)}
{\partial x_t(s,a)}
=
\frac{\mathbf{1}_{\{a=b\}}-\pi_t(s,b)}
{\mu_t(s)}.
\end{align}
Therefore, differentiating the stage potential yields
\begin{align}
\frac{\partial\phi(x_t)}{\partial x_t(s,a)}
&=r(s,a,x_t)- H_s^*(\pi_t(s,\cdot))
- \sum_{a'} \frac{\partial H_s^*(\pi_t(s,\cdot))}{\partial\pi_t(s,a')}\left(\mathbf{1}_{\{a'=a\}}-\pi_t(s,a')\right)\nonumber\\
&=r(s,a,x_t)-H_s^*(\pi_t(s,\cdot)) 
-\frac{\partial H_s^*(\pi_t(s,\cdot))}{\partial\pi_t(s,a)}+\pi_t(s,\cdot)^\top\nabla H_s^*(\pi_t(s,\cdot))\nonumber\\
&=R(s,a,x_t).
\end{align}
At $\mu_t(s)=0$, the corresponding term
$\mu_t(s)H_s^*(\pi_t(s,\cdot))$ is defined
by continuity as zero.
This establishes the stated potential
representation at states with positive
population mass.

\section{Structural property of transition kernel}\label{app:trans_kernel}

We first show that the transition kernel proposed in Eq.~\eqref{eq:num:transition} is valid with $\kappa_1, \kappa_2\geq 0, \kappa_1+\kappa_2< 1$. 

The factor of the first term in Eq.~\eqref{eq:num:action-component} is bounded as $-1\leq 2y_{t,a}/\|y_t\|_2-1\leq1$, while the double summation in the second term in Eq.~\eqref{eq:num:action-component} is bounded as follows:
\begin{align}
0\leq\sum_{1\leq i<j\leq m}(y_{t,i}-y_{t,j})^2
&=\frac12\sum_{i=1}^{m}\sum_{j=1}^{m}
\left(y_{t,i}^2+y_{t,j}^2-2y_{t,i}y_{t,j}\right)=m\sum_{a=1}^{m}y_{t,a}^2
-\left(\sum_{a=1}^{m}y_{t,a}\right)^2\nonumber\\
&=m\sum_{a=1}^{m}y_{t,a}^2-1\leq m-1,
\end{align}
where the third equality is due to $\sum_a y_{t,a}=1$ and the inequality is due to $\sum_a y_{t,a}^2\leq (\sum_a y_{t,a})^2 = 1$. 

Accordingly, we have $-\kappa_1\leq q_a(y_t) \leq \kappa_1+\kappa_2$. The condition of $\kappa_1, \kappa_2$ thus naturally ensures $0< P(s'|s,a,x_t)< 1$ for any $s',s,a,x_t$.

To validate Assumption~\ref{ass:transition}, we first differentiate the two terms in Eq.~\eqref{eq:num:action-component} with respect to $y_{t,a'}$ for any $a'\in\setA$, which yields
\begin{align}
&\frac{\partial}{\partial y_{t,a'}}
\left(\frac{2y_{t,a}}{\|y_t\|}-1\right)=
2\left(\frac{\mathbf 1_{\{a=a'\}}}{\|y_t\|}-
\frac{y_{t,a}y_{t,a'}}{\|y_t\|^3}
\right),\label{eq:num:partial_1} \\
&\frac{\partial}{\partial y_{t, a'}}\left[\frac{\sum_{a_1=1}^{m-1}\sum_{a_2=a_1+1}^m (y_{t,a_1}-y_{t,a_2})^2}{m-1}\right]=\frac{2}{m-1}(my_{t,a'}-1).\label{eq:num:partial_2}
\end{align}

Plugging Eqs.~\eqref{eq:num:partial_1} and \eqref{eq:num:partial_2} into Eq.~\eqref{eq:expect_trans_marg} yields
\begin{align}
    \mathbb{E}\left[\frac{\partial P(s'|\cdot, x_t)}{\partial x_t(s,a)}\right] &= \sum_{s''\in\setS}\sum_{a''\in\setA} \frac{\partial P(s'|s'', a'', x_t)}{\partial x_t(s,a)}x_t(s'',a'') = \sum_{s''\in\setS}\sum_{a''\in\setA} \frac{\partial P(s'|s'', a'', x_t)}{\partial y_{t,a}}\frac{\partial y_{t,a}}{\partial x_t(s,a)} x_t(s'', a'')\nonumber\\
    &= \frac{z_{s'}}{n}\sum_{a''\in\setA} \left[2\kappa_1 \left(\frac{\mathbf 1_{\{a''=a\}}}{\|y_t\|}-\frac{y_{t,a''}y_{t,a}}{\|y_t\|^3}\right) + \frac{2\kappa_2}{m-1}(my_{t,a}-1) \right] \sum_{s''\in\setS} x_t(s'',a'')\nonumber\\
    &=\frac{2\kappa_1z_{s'}}{n} \sum_{a''\in\setA} \left(\frac{\mathbf 1_{\{a''=a\}}}{\|y_t\|}-\frac{y_{t,a''}y_{t,a}}{\|y_t\|^3}\right) y_{t,a''} + \frac{2\kappa_2z_{s'}}{n(m-1)} (my_{t,a}-1)\sum_{a''\in\setA}y_{t,a''}\nonumber\\
    &= \frac{2\kappa_2z_{s'}}{n(m-1)} (my_{t,a}-1)
\end{align}
Therefore, Assumption~\ref{ass:transition} is violated as long as $\kappa_2>0$.

\section{Other numerical results}\label{app:numerical}
\begin{figure}
    \centering
    \includegraphics[width=\linewidth]{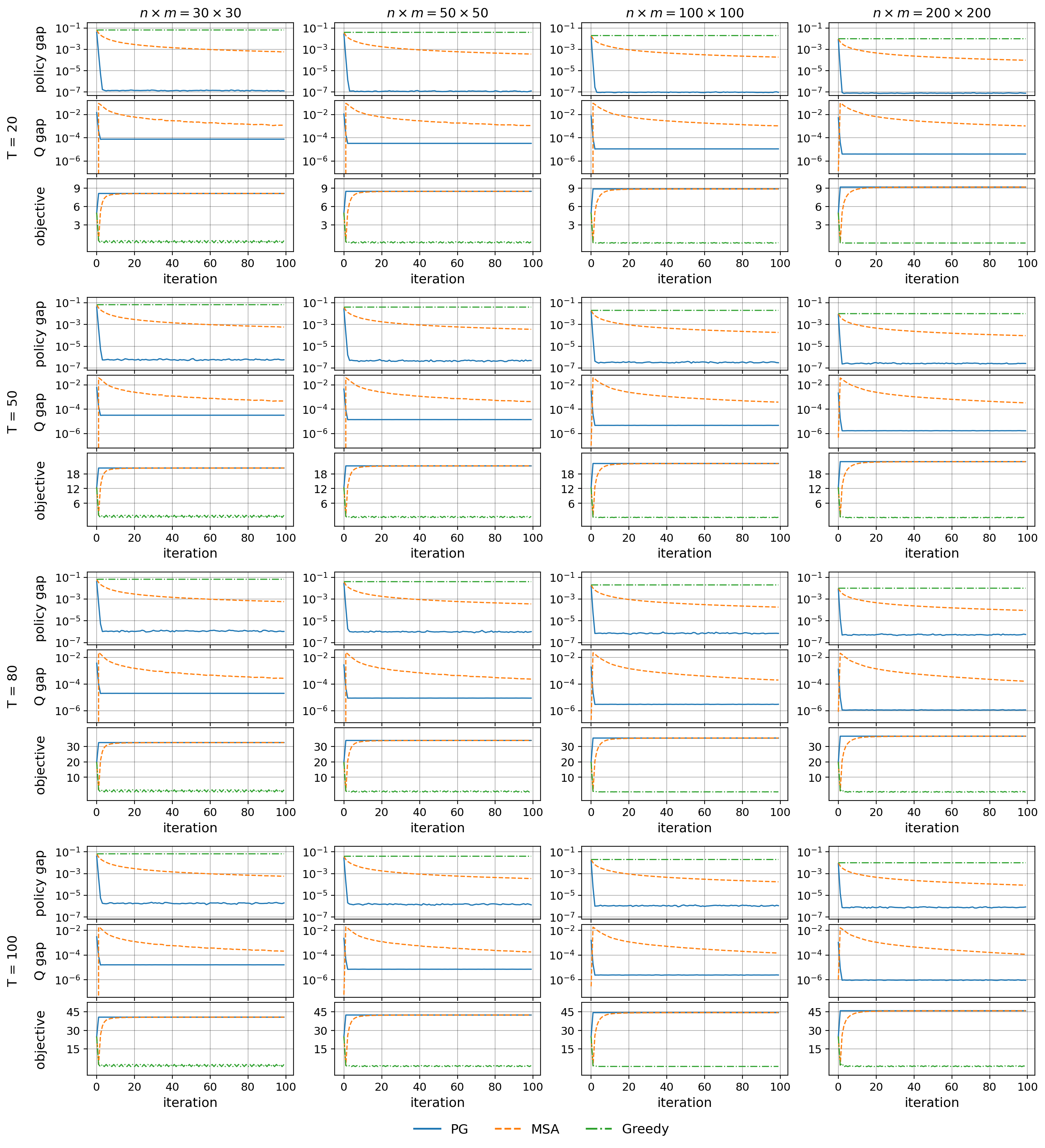}
    \caption{All convergence results}
    \label{fig:conv_scale_all}
\end{figure}

Figure~\ref{fig:conv_scale_all} compares the convergence of PG, MSA, and Greedy methods across different values of
$n\times m$ and $T$. Throughout the tested range,
PG consistently achieves small policy and $Q$-gaps in
substantially fewer iterations.

\end{document}